\documentclass[twoside,12pt]{article}
\usepackage{amssymb}
\usepackage{a4wide}
\usepackage{amsthm}
\usepackage{amsmath}
\usepackage[all]{xy}
\usepackage{enumerate}
\usepackage{fancyhdr}

\newtheorem{theorem}{Theorem}[section]
\newtheorem{proposition}[theorem]{Proposition}
\newtheorem{corollary}[theorem]{Corollary}
\newtheorem{lemma}[theorem]{Lemma}
\theoremstyle{definition}
\newtheorem*{notation}{Notation}
\newtheorem*{Beweis}{Proof}
\newtheorem{definition}[theorem]{Definition}
\newtheorem{punto}[theorem]{}
\theoremstyle{remark}
\newtheorem{remark}[theorem]{Remark}
\newtheorem{ex}[theorem]{Example}
\newtheorem{remarks}[theorem]{Remarks}
\swapnumbers
\input{tcilatex}

\begin{document}

\title{\textbf{Injective Morita Contexts (Revisited)}\bigskip \\
{\large {\emph{Dedicated to Prof. Robert Wisbauer}}}}
\author{\begin{tabular}{lll}
{\bf J. Y. Abuhlail} \thanks{
Corresponding Author} &  & {\bf S. K. Nauman} \\
Department of Mathematics $\&$ Statistics &  & Department of Mathematics \\
King Fahd University of Petroleum &  & King AbdulAziz University \\
\&\ Minerals, Box $\#$\ 5046 &  & P.O.Box 80203 \\
31261 Dhahran (KSA) &  & 21589 Jeddah (KSA) \\
abuhlail@kfupm.edu.sa &  & synakhaled@hotmail.com
\end{tabular}
}

\date{}
\maketitle

\begin{abstract}
This paper is an exposition of the so-called \emph{injective Morita contexts}
(in which the connecting bimodule morphisms are injective) and \emph{Morita }%
$\alpha $\emph{-contexts} (in which the connecting bimodules enjoy some
local projectivity in the sense of Zimmermann-Huisgen). Motivated by
situations in which only one trace ideal is in action, or the compatibility
between the bimodule morphisms is not needed, we introduce the notions of
Morita \emph{semi-contexts} and \emph{Morita data}, and investigate them.
Injective Morita data will be used (with the help of \emph{static }and \emph{%
adstatic modules})\emph{\ }to establish equivalences between some \emph{%
intersecting subcategories} related to subcategories of categories of
modules that are localized or colocalized by trace ideals of a Morita datum.
We end up with applications of Morita $\alpha $-contexts to $\ast $\emph{%
-modules} and \emph{injective right wide Morita contexts.}
\end{abstract}

\section{Introduction}

\qquad \emph{Morita contexts}, in general, and (\emph{semi-})\emph{strict
Morita} contexts (with surjective connecting bilinear morphisms), in
particular, were extensively studied and developed exponentially during the
last few decades (e.g. \cite{AGH-Z1997}). However, we sincerely feel that
there is a gap in the literature on \emph{injective Morita contexts} (i.e.
those with injective connecting bilinear morphisms). Apart from the results
in \cite{Nau-1994-a}, \cite{Nau-1994-b} (where the second author initially
explored this notion) and from an application to Grothendieck groups in the
recent paper (\cite{Nau2004}), it seems that injective Morita contexts were
not studied \emph{systematically} at all.

We noticed that in several results of (\cite{Nau-1993}, \cite{Nau-1994-a}
and \cite{Nau-1994-b}) that are related to Morita contexts, only one trace
ideal is used. Observing this fact, we introduce the notions of \emph{Morita
semi-contexts} and \emph{Morita data} and investigate them. Several results
are proved then for \emph{injective} Morita semi contexts and/or injective
Morita data.

Consider a Morita datum $\mathcal{M}=(T,S,P,Q,<,>_{T},<,>_{S}),$ with not
necessarily compatible bimodule morphisms $<,>_{T}:P\otimes _{S}Q\rightarrow
T$ and $<,>_{S}:Q\otimes _{T}P\rightarrow S.$ We say that $\mathcal{M}$ is
\emph{injective}, iff $<,>_{T}$ and $<,>_{S}$ are injective, and to be a
\emph{Morita }$\alpha $\emph{-datum}, iff the associated dual pairings $%
\mathbf{P}_{l}:=(Q,$ $_{T}P),$ $\mathbf{P}_{r}:=(Q,P_{S}),$ $\mathbf{Q}%
_{l}:=(P,$ $_{S}Q)$ and $\mathbf{Q}_{r}:=(P,Q_{T})$ satisfy the $\alpha $%
-condition (which is closely related to the notion of local projectivity in
the sense of Zimmermann-Huisgen \cite{Z-H-1976}). The $\alpha $-condition
was introduced in \cite{AG-TL2001} and further investigated by the first
author in \cite{Abu-2005}.

While (semi-)strict unital Morita contexts induce equivalences between the
whole module categories of the rings under consideration, we show in this
paper how injective Morita (semi-)contexts and injective Morita data play an
important role in establishing equivalences between suitable \emph{%
intersecting subcategories} of module categories (e.g. intersections of
subcategories that are localized/colocalized by trace ideals of a Morita
datum with subcategories of static/adstatic modules, etc.). Our main
applications in addition to equivalences related to the Kato-Ohtake-M\"{u}%
ller \emph{localization-colocalization theory} (developed in \cite{Kat-1978}%
, \cite{KO-1979} and \cite{Mul-1974}), will be to $\ast $\emph{-modules}
(introduced by Menini and Orsatti \cite{MO-1989}) and to \emph{right} \emph{%
wide Morita contexts} (introduced by F. Casta\~{n}o Iglesias and J. G\'{o}%
mez-Torrecillas \cite{C-IG-T1995}).

Most of our results will be stated for \emph{left modules}, while deriving
the \textquotedblleft dual\textquotedblright\ versions for right modules is
left to the interested reader. Moreover, for Morita contexts, some results
are stated/proved for only one of the Morita semi-contexts, as the ones
corresponding to the second semi-context can be obtained analogously. For
the convenience of the reader, we tried to make the paper self-contained, so
that it can serve as a reference on injective \emph{Morita }(\emph{semi-})%
\emph{contexts} and their applications. In this respect, and for the sake of
completeness, we have included some previous results of the authors that are
(in most cases) either provided with new shorter proofs, or are obtained
under weaker conditions.

This paper is organized as follows: After this brief introduction, we give
in Section 2 some preliminaries including the basic properties of \emph{dual
}$\alpha $\emph{-pairings,} which play a central role in rest of the work.
The notions of \emph{Morita semi-contexts} and \emph{Morita data} are
introduced in Section 3, where we clarify their relations with the \emph{%
dual pairings} and the so-called \emph{elementary rngs}. \emph{Injective
Morita }(\emph{semi-})\emph{contexts} appear in Section 4, where we study
their interplay with dual $\alpha $-pairings and provide some examples and a
counter-example. In Section 5 we include some observations regarding \emph{%
static} and \emph{adstatic} modules and use them to obtain equivalences
among suitable \emph{intersecting subcategories} of modules related to a
Morita (semi-)context. In the last section, more applications are presented,
mainly to subcategories of modules that are \emph{localized} or \emph{%
colocalized} by a trace ideal of an injective Morita (semi-)context, to $%
\ast $\emph{-modules} and to \emph{injective right wide Morita contexts}.

\section{Preliminaries}

\qquad Throughout, $R$ denotes a commutative ring with $1_{R}\neq 0_{R}$ and
$A,A^{\prime },B,B^{\prime }$ are unital $R$-algebras. We have reserved the
term \textquotedblleft \emph{ring}\textquotedblright\ for an associative
ring with a multiplicative unity, and we will use the term \textquotedblleft
\emph{rng}\textquotedblright\ for a general associative ring (not
necessarily with unity). All modules over rings are assumed to be unitary,
and ring morphisms are assumed to respect multiplicative unities. If $%
\mathfrak{T}$ and $\mathfrak{S}$ are categories, then we write $\mathfrak{T}%
\leq \mathfrak{S}$ ($\mathfrak{T}\leq \mathfrak{S}$) to mean that $\mathfrak{%
T}$ is a (full) subcategory of $\mathfrak{S},$ and $\mathfrak{T}\approx
\mathfrak{S}$ to indicate that $\mathfrak{T}$ and $\mathfrak{S}$ are
equivalent.

\subsection*{Rngs and their modules}

\begin{punto}
By an $A$\textbf{-rng} $(T,\mu _{T}),$ we mean an $(A,A)$-bimodule $T$ with
an $(A,A)$-bilinear morphism $\mu _{T}:T\otimes _{A}T\rightarrow T,$ such
that $\mu _{T}\circ (\mu _{T}\otimes _{A}id_{T})=\mu _{T}\circ
(id_{T}\otimes _{A}\mu _{T}).$ We call an $A$-rng $(T,\mu _{T})$ an $A$%
\textbf{-ring}, iff there exists in addition an $(A,A)$-bilinear morphism $%
\eta _{T}:A\rightarrow T,$ called the \textbf{unity map}, such that $\mu
_{T}\circ (\eta _{T}\otimes _{A}id_{T})=\vartheta _{T}^{l}$ and $\mu
_{T}\circ (id_{T}\otimes _{A}\eta _{T})=\vartheta _{T}^{r}$ (where $A\otimes
_{A}T\overset{\vartheta _{T}^{l}}{\simeq }T$ and $T\otimes _{A}A\overset{%
\vartheta _{T}^{r}}{\simeq }T$ are the canonical isomorphisms). So, an $A$%
-ring is a unital $A$-rng; and an $A$-rng is (roughly speaking) an $A$-ring
not necessarily with unity.
\end{punto}

\begin{punto}
A morphism of rngs $(\psi :\delta ):(T:A)\rightarrow (T^{\prime }:A^{\prime
})$ consists of a morphism of $R$-algebras $\delta :A\rightarrow A^{\prime }$
and an $(A,A)$-bilinear morphism $\psi :T\rightarrow T^{\prime },$ such that
$\mu _{T^{\prime }}\circ \chi _{(T^{\prime },T^{\prime })}^{(A,A^{\prime
})}\circ (\psi \otimes _{A}\psi )=\psi \circ \mu _{T}$ (where $\chi
_{(T^{\prime },T^{\prime })}^{(A,A^{\prime })}:T^{\prime }\otimes
_{A}T^{\prime }\rightarrow T^{\prime }\otimes _{A^{\prime }}T^{\prime }$ is
the canonical map induced by $\delta $). By $\mathbb{RNG}$ we denote the
category of associative rngs with morphisms being rng morphisms, and by $%
\mathbb{URNG}<\mathbb{RNG}$ the (non-full) subcategory of \emph{unital}
rings with morphisms being the morphisms in $\mathbb{RNG}$ which respect
multiplicative unities.
\end{punto}

\begin{punto}
Let $(T,\mu _{T})$ be an $A$-rng. By a \textbf{left }$T$\textbf{-module} we
mean a left $A$-module $N$ with a left $A$-linear morphism $\phi
_{T}^{N}:T\otimes _{A}N\rightarrow N,$ such that $\phi _{T}^{N}\circ (\mu
_{T}\otimes _{A}id_{N})=\phi _{T}^{N}\circ (id_{T}\otimes _{A}\phi
_{T}^{N}). $ For left $T$-modules $M,N,$ we call a left $A$-linear morphism $%
f:M\rightarrow N$ a $T$\textbf{-linear morphism}, iff $f(tm)=tf(m)$ for all $%
t\in T.$ The category of left $T$-modules and left $T$-linear morphisms is
denoted by $_{T}\mathbb{M}.$ The category $\mathbb{M}_{T}$ of right $T$%
-modules is defined analogously. Let $(T:A)$ and $(T^{\prime }:A^{\prime })$
be rngs. We call an $(A,A^{\prime })$-bimodule $N$ a $(T,T^{\prime })$%
\textbf{-bimodule}, iff $(N,\phi _{T}^{N})$ is a left $T$-module and $%
(N,\phi _{T^{\prime }}^{N})$ is a right $T^{\prime }$-module, such that $%
\phi _{T^{\prime }}^{N}\circ (\phi _{T}^{N}\otimes _{A^{\prime
}}id_{T^{\prime }})=\phi _{T}^{N}\circ (id_{T}\otimes _{A}\phi _{T^{\prime
}}^{N}).$ For $(T,T^{\prime })$-bimodules $M,N,$ we call an $(A,A^{\prime })$%
-bilinear morphism $f:M\rightarrow N$ $(T,T^{\prime })$\textbf{-bilinear},
provided $f$ is left $T$-linear and right $T^{\prime }$-linear. The category
of $(T,T^{\prime })$-bimodules is denoted by $_{T}\mathbb{M}_{T^{\prime }}.$
In particular, for any $A$-rng $T,$ a left (right) $T$-module $M$ has a
canonical structure of a \emph{unitary} right (left) $S$-module, where $S:=%
\mathrm{End}(_{T}M)^{op}$ ($S:=\mathrm{End}(M_{T})$); and moreover, with
this structure $M$ becomes a $(T,S)$-bimodule (an $(S,T)$-bimodule).
\end{punto}

\begin{remark}
Similarly, one can define rngs over arbitrary (not-necessarily unital)
ground rngs and rng morphisms between them. Moreover, one can define
(bi)modules over such rngs and (bi)linear morphisms between them.
\end{remark}

\begin{notation}
Let $T$ be an $A$-rng. We write $_{T}U$ ($U_{T}$) to denote that $U$ is a
left (right) $T$-module. For a left (right) $T$-module $_{T}U,$ we consider
the set $^{\ast }U:=\mathrm{Hom}_{T-}(U,T)$ ($U^{\ast }:=\mathrm{Hom}%
_{-T}(U,T)$) of all left (right) $T$-linear morphisms from $U$ to $T$ with
the canonical right (left) $T$-module structure.
\end{notation}

\subsection*{Generators and cogenerators}

\begin{definition}
Let $T$ be an $A$-rng. For a left $T$-module $_{T}U$ consider the following
subclasses of $_{T}\mathbb{M}:$%
\begin{equation*}
\begin{tabular}{lll}
$\mathrm{Gen}(_{T}U)$ & $:=$ & $\{_{T}V\mid \exists $ a set $\Lambda $ and
an exact sequence $U^{(\Lambda )}\rightarrow V\rightarrow 0\};$ \\
$\mathrm{Cogen}(_{T}U)$ & $:=$ & $\{_{T}W\mid \exists $ a set $\Lambda $ and
an exact sequence $0\rightarrow W\rightarrow U^{\Lambda }\};$ \\
$\mathrm{Pres}(_{T}U)$ & $:=$ & $\{_{T}V\mid \exists $ sets $\Lambda
_{1},\Lambda _{2}$ and an exact sequence $U^{(\Lambda _{2})}\rightarrow
U^{(\Lambda _{1})}\rightarrow V\rightarrow 0\};$ \\
$\mathrm{Copres}(_{T}U)$ & $:=$ & $\{_{T}W\mid \exists $ sets $\Lambda
_{1},\Lambda _{2}$ and an exact sequence $0\rightarrow W\rightarrow
U^{\Lambda _{1}}\rightarrow U^{\Lambda _{2}}\};$%
\end{tabular}%
\end{equation*}%
A left $T$-module in $\mathrm{Gen}(_{T}U)$ (respectively $\mathrm{Cogen}%
(_{T}U),$ $\mathrm{Pres}(_{T}U),$ $\mathrm{Copres}(_{T}U)$) is said to be $U$%
\textbf{-generated }(respectively $U$\textbf{-cogenerated}, $U$\textbf{%
-presented}, $U$\textbf{-copresented}). Moreover, we say that $_{T}U$ is a
\textbf{generator} (respectively \textbf{cogenerator}, \textbf{presentor},
\textbf{copresentor}), iff $\mathrm{Gen}(_{T}U)=$ $_{T}\mathbb{M}$
(respectively $\mathrm{Cogen}(_{T}U)=$ $_{T}\mathbb{M},$ $\mathrm{Pres}%
(_{T}U)=$ $_{T}\mathbb{M},$ $\mathrm{Copres}(_{T}U)=$ $_{T}\mathbb{M}$).
\end{definition}

\qquad

\subsection*{Dual $\protect\alpha $-pairings}

\qquad In what follows we recall the definition and properties of dual $%
\alpha $-pairings introduced in \cite[Definition 2.3.]{AG-TL2001} and
studied further in \cite{Abu-2005}.

\begin{punto}
Let $T$ be an $A$-rng. A \textbf{dual left }$T$\textbf{-pairing} $\mathbf{P}%
_{l}=(V,$ $_{T}W)$ consists of a left $T$-module $W$ and a right $T$-module $%
V$ with a right $T$-linear morphism $\kappa _{\mathbf{P}_{l}}:V\rightarrow $
$^{\ast }W$ (equivalently a left $T$-linear morphism $\chi _{\mathbf{P}%
_{l}}:W\rightarrow V^{\ast }$). For dual left pairings $\mathbf{P}_{l}=(V,$ $%
_{T}W),$ $\mathbf{P}_{l}^{\prime }=(V^{\prime },$ $_{T^{\prime }}W^{\prime
}),$ a morphism of dual left pairings $(\xi ,\theta ):(V^{\prime },W^{\prime
})\rightarrow (V,W)$ consists of a triple%
\begin{equation*}
(\xi ,\theta :\varsigma ):(V,\text{ }_{T}W)\rightarrow (V^{\prime },\text{ }%
_{T^{\prime }}W^{\prime }),
\end{equation*}%
where $\xi :V\rightarrow V^{\prime }$ and $\theta :W^{\prime }\rightarrow W$
are $T$-linear and $\varsigma :T\rightarrow T^{\prime }$ is a morphism of
rngs, such that considering the induced maps $<,>_{T}:V\times W\rightarrow T$
and $<,>_{T^{\prime }}:V^{\prime }\times W^{\prime }\rightarrow T^{\prime }$
we have%
\begin{equation}
<\xi (v),w^{\prime }>_{T^{\prime }}=\varsigma (<v,\theta (w^{\prime })>_{T})%
\text{ for all }v\in V\text{ and }w^{\prime }\in W^{\prime }.  \label{comp}
\end{equation}%
The dual left pairings with the morphisms defined above build a category,
which we denote by $\mathcal{P}_{l}.$ With $\mathcal{P}_{l}(T)\leq \mathcal{P%
}_{l}$ we denote the full subcategory of dual $T$-pairings. The category $%
\mathcal{P}_{r}$ of dual right pairings and its full subcategory $\mathcal{P}%
_{r}(T)\leq \mathcal{P}_{r}$ of dual right $T$-pairings are defined
analogously.
\end{punto}

\begin{remark}
The reader should be warned that (in general) for a non-commutative rng $T$
and a dual left $T$-pairing $\mathbf{P}_{l}=(V,$ $_{T}W),$ the following map
induced by the right $T$-linear morphism $\kappa _{\mathbf{P}%
_{l}}:V\rightarrow $ $^{\ast }W:$%
\begin{equation*}
<,>_{T}:V\times W\rightarrow T,\text{ }<v,w>_{T}:=\kappa _{\mathbf{P}%
_{l}}(v)(w)
\end{equation*}%
is not necessarily $T$\emph{-balanced}, and so does not induce (in general)
a map $V\otimes _{T}W\rightarrow T.$ In fact, for all $v\in V,$ $w\in W$ and
$t\in T$ we have%
\begin{equation*}
\begin{tabular}{lllllllll}
$<vt,w>$ & $=$ & $\kappa _{\mathcal{P}_{l}}(vt)(w)$ & $=$ & $[\kappa _{%
\mathcal{P}_{l}}(v)t](w)$ & $=$ & $[\kappa _{\mathcal{P}_{l}}(v)(w)]t$ & $=$
& $<v,w>_{T}t;$ \\
$<v,tw>$ & $=$ & $\kappa _{\mathcal{P}_{l}}(v)(tw)$ & $=$ & $t[\kappa _{%
\mathcal{P}_{l}}(v)(w)]$ & $=$ & $t<v,w>_{T}.$ &  &
\end{tabular}%
\end{equation*}
\end{remark}

\begin{punto}
Let $T$ be an $A$-rng, $N,W$ be left $T$-modules and identify $N^{W}$ with
the set of all mappings from $W$ to $N.$ Considering $N$ with the \emph{%
discrete topology} and $N^{W}$ with the product topology, the induced \emph{%
relative topology }on $\mathrm{Hom}_{T-}(W,N)\hookrightarrow N^{W}$ is a
linear topology (called the \textbf{finite topology}), for which the \emph{%
basis of neighborhoods of }$0$ is given by the set of annihilator submodules:%
\begin{equation*}
\mathcal{B}_{f}{\normalsize (0)}:=\{F^{\bot (\mathrm{Hom}_{T-}(W,N))}\mid
F=\{w_{1},...,w_{k}\}\subset W\text{ is a finite subset}\},
\end{equation*}%
where%
\begin{equation*}
F^{\bot (\mathrm{Hom}_{T-}(W,N))}:=\{f\in \mathrm{Hom}_{T-}(W,N))\mid
f(W)=0\}.
\end{equation*}
\end{punto}

\begin{punto}
\label{P-alph}Let $T$ be an $A$-rng, $\mathbf{P}_{l}=(V,$ $_{T}W)$ a dual
left $T$-pairing and consider for every right $T$-module $U_{T}$ the
following canonical map%
\begin{equation}
\alpha _{U}^{\mathbf{P}_{l}}:U\otimes _{T}W\rightarrow \mathrm{Hom}%
_{-T}(V,U),\text{ }\sum u_{i}\otimes _{T}w_{i}\mapsto \lbrack v\mapsto \sum
u_{i}<v,w_{i}>_{T}].  \label{alp}
\end{equation}%
We say that $\mathbf{P}_{l}=(V,$ $_{T}W)\in \mathcal{P}_{l}(T)$ \textbf{%
satisfies the left }$\alpha $\textbf{-condition} (or is a \textbf{dual left }%
$\alpha $\textbf{-pairing}), iff $\alpha _{U}^{\mathbf{P}_{l}}$ is injective
for every right $T$-module $U_{T}.$ By $\mathcal{P}_{l}^{\alpha }(T)\leq
\mathcal{P}_{l}(T)$ we denote the \emph{full} subcategory of dual left $T$%
-pairings satisfying the left $\alpha $-condition. The full subcategory of
\textbf{dual right }$\alpha $\textbf{-pairings} $\mathcal{P}_{r}^{\alpha
}(T)\leq \mathcal{P}_{r}(T)$ is defined analogously.
\end{punto}

\begin{definition}
Let $T$ be an $A$-rng, $\mathbf{P}_{l}=(V,$ $_{T}W)$ be a dual left $T$%
-pairing and consider%
\begin{equation*}
\kappa _{\mathbf{P}_{l}}:V\rightarrow \text{ }^{\ast }W\text{ and }\alpha
_{V}^{\mathbf{P}_{l}}:V\otimes _{T}W\rightarrow \mathrm{End}(V_{T}).
\end{equation*}%
We say $\mathbf{P}_{l}\in \mathcal{P}_{l}(T)$ is

\textbf{dense,} iff $\kappa _{\mathbf{P}_{l}}(V)\subseteq $ $^{\ast }W$ is
dense (w.r.t. the \emph{finite topology} on $^{\ast }W\hookrightarrow T^{W}$%
\emph{)};

\textbf{injective} (resp. \textbf{semi-strict}, \textbf{strict}), iff $%
\alpha _{V}^{\mathbf{P}_{l}}$ is injective (resp. surjective, bijective);

\textbf{non-degenerate}, iff $V\overset{\kappa _{\mathbf{P}_{l}}}{%
\hookrightarrow }$ $^{\ast }W$ and $W\overset{\chi _{\mathbf{P}_{l}}}{%
\hookrightarrow }V^{\ast }$ canonically.
\end{definition}

\begin{punto}
Let $T$ be an $A$-rng. We call a $T$-module $W$ \textbf{locally projective}
(in the sense of B. Zimmermann-Huisgen \cite{Z-H-1976}), iff for every
diagram of $T$-modules%
\begin{equation*}
\xymatrix{0 \ar[r] & F \ar@{.>}[dr]_{g' \circ \iota} \ar[r]^{\iota} & W
\ar[dr]^{g} \ar@{.>}[d]^{g'} & & \\ & & L \ar[r]_{\pi} & N \ar[r] & 0}
\end{equation*}%
with exact rows and finitely generated $T$-submodule $F\subseteq W$ : for
every $T$-linear morphism $g:W\rightarrow N,$ there exists a $T$-linear
morphism $g^{\prime }:W\rightarrow L,$ such that $g\circ \iota =\pi \circ
g^{\prime }\circ \iota .$
\end{punto}

\qquad For proofs of the following basic properties of \emph{locally
projective modules} and \emph{dual }$\alpha $\emph{-pairings }see \cite%
{Abu-2005} and \cite{Z-H-1976}:

\begin{proposition}
\label{rp-rp}\emph{\ }Let $T$ be an $A$-ring and $\mathbf{P}_{l}=(V,$ $%
_{T}W)\in \mathcal{P}_{l}(T).$

\begin{enumerate}
\item The left $T$-module $_{T}W$ is locally projective if and only if $%
(^{\ast }W,W)$ is an $\alpha $-pairing.

\item The left $T$-module $_{T}W$ is locally projective, iff for any finite
subset $\{w_{1},...,w_{k}\}\subseteq W,$ there exists $\{(f_{i},\widetilde{w}%
_{i})\}_{i=1}^{k}\subset $ $^{\ast }W\times W$ such that $%
w_{j}=\dsum\limits_{i=1}^{k}f_{i}(w_{j})\widetilde{w}_{i}$ for all $%
j=1,...,k.$

\item If $_{T}W$ is locally projective, then $_{T}W$ is flat and $T$%
-cogenerated.

\item If $\mathbf{P}_{l}\in \mathcal{P}_{l}^{\alpha }(T),$ then $_{T}W$ is
locally projective.

\item If $_{T}W$ is locally projective and $\kappa _{P}(V)\subseteq $ $%
^{\ast }W$ is dense, then $\mathbf{P}_{l}\in \mathcal{P}_{l}^{\alpha }(T).$

\item Assume $T_{T}$ is an injective cogenerator. Then $\mathbf{P}_{l}\in
\mathcal{P}_{l}^{\alpha }(T)$ if and only if $_{T}W$ is locally projective
and $\kappa _{\mathbf{P}_{l}}(V)\subseteq $ $^{\ast }W$ is dense.

\item If $T$ is a QF ring, then $\mathbf{P}_{l}\in \mathcal{P}_{l}^{\alpha
}(T)$ if and only if $_{T}W$ is projective and $W\overset{\chi _{\mathbf{P}%
_{l}}}{\hookrightarrow }V^{\ast }.$
\end{enumerate}
\end{proposition}

\qquad The following result completes the nice observation \cite[42.13.]%
{BW-2003} about locally projective modules:

\begin{proposition}
\label{lp-dense}Let $T$ be a ring, $_{T}W$ a left $T$-module, $S:=\mathrm{End%
}(_{T}W)^{op}$ and consider the canonical $(S,S)$-bilinear morphism%
\begin{equation*}
\lbrack ,]_{W}:\text{ }^{\ast }W\otimes _{T}W\rightarrow \mathrm{End}(_{T}W),%
\text{ }f\otimes _{T}w\mapsto \lbrack \widetilde{w}\mapsto f(\widetilde{w}%
)w].
\end{equation*}

\begin{enumerate}
\item $_{T}W$ is \emph{finitely generated projective} if and only if $%
[,]_{W} $ is surjective.

\item $_{T}W$ is locally projective if and only if $\mathrm{Im}%
([,]_{W})\subseteq \mathrm{End}(_{T}W)$ is dense.
\end{enumerate}
\end{proposition}

\begin{Beweis}
\begin{enumerate}
\item This follows by \cite[12.8.]{Fai-1981}.

\item Assume $_{T}W$ is locally projective and consider for every left $T$%
-module $N$ the canonical mapping
\begin{equation*}
\lbrack ,]_{N}^{W}:^{\ast }W\otimes _{T}N\rightarrow \mathrm{Hom}_{T}(W,N),%
\text{ }f\otimes _{T}n\mapsto \lbrack \widetilde{w}\mapsto f(\widetilde{w}%
)n].
\end{equation*}%
It follows then by \cite[42.13.]{BW-2003}, that $\mathrm{Im}%
([,]_{N}^{W})\subseteq \mathrm{Hom}_{T}(W,N)$ is dense. In particular,
setting $N=W$ we conclude that $\mathrm{Im}([,]_{W})\subseteq \mathrm{End}%
(_{T}W)$ is dense. On the other hand, assume $\mathrm{Im}([,]_{W})\subseteq
\mathrm{End}(_{T}W)$ is dense. Then for every finite subset $%
\{w_{1},...,w_{k}\}\subseteq W,$ there exists $\dsum\limits_{i=1}^{n}%
\widetilde{g}_{i}\otimes _{T}\widetilde{w}_{i}\in $ $^{\ast }W\otimes _{T}W$
with%
\begin{equation*}
w_{j}=id_{W}(w_{j})=[,]_{W}(\dsum_{i=1}^{n}\widetilde{g}_{i}\otimes _{T}%
\widetilde{w}_{i})(w_{j})=\dsum_{i=1}^{n}\widetilde{g}_{i}(w_{j})\widetilde{w%
}_{i}\text{ for }j=1,...,k.
\end{equation*}%
It follows then by Proposition \ref{rp-rp} \textquotedblleft
2\textquotedblright\ that $_{T}W$ is locally projective.$\blacksquare $
\end{enumerate}
\end{Beweis}

\section{Morita (Semi)contexts}

\qquad We noticed, in the proofs of some results on equivalences between
subcategories of module categories associated to a given Morita context,
that no use is made of the \emph{compatibility} between the connecting
bimodule morphisms (or even that only one trace ideal is used and so only
one of the two bilinear morphisms is really in action). Some results of this
type appeared, for example, in \cite{Nau-1993}, \cite{Nau-1994-a} and \cite%
{Nau-1994-b}. Moreover, in our considerations some Morita contexts will be
formed for arbitrary associative rngs (i.e. not necessarily unital rings).
These considerations motivate us to make the following general definitions:

\begin{punto}
By a \textbf{Morita semi-context }we mean a tuple%
\begin{equation}
\mathbf{m}_{T}=((T:A),(S:B),P,Q,<,>_{T},I),  \label{semi-Mor}
\end{equation}%
where $T$ is an $A$-rng, $S$ is a $B$-rng, $P$ is a $(T,S)$-bimodule, $Q$ is
an $(S,T)$-bimodule, $<,>_{T}:P\otimes _{S}Q\rightarrow T$ is a $(T,T)$%
-bilinear morphism and $I:=\func{Im}(<,>_{T})\vartriangleleft T$ (called the
\textbf{trace ideal} \textbf{associated to} $\mathbf{m}_{T}).$ We drop the
ground rings $A,B$ and the trace ideal $I\vartriangleleft T,$ if they are
not explicitly in action. If $\mathbf{m}_{T}$ (\ref{semi-Mor}) is a Morita
semi-context and $T,S$ are unital rings, then we call $\mathbf{m}_{T}$ a
\textbf{unital Morita semi-context}.
\end{punto}

\begin{punto}
Let $\mathbf{m}_{T}=((T:A),(S:B),P,Q,<,>_{T}),$ $\mathbf{m}_{T^{\prime
}}=((T^{\prime }:A^{\prime }),(S^{\prime }:B^{\prime }),P^{\prime
},Q^{\prime },<,>_{T^{\prime }})$ be Morita semi-contexts. By a \textbf{%
morphism of Morita semi-contexts} from $\mathbf{m}_{T}$ to $\mathbf{m}%
_{T^{\prime }}$ we mean a four fold set of morphisms
\begin{equation*}
((\beta :\delta ),(\gamma :\sigma ),\phi ,\psi
):((T:A),(S:B),P,Q)\rightarrow ((T^{\prime }:A^{\prime }),(S^{\prime
}:B^{\prime }),P^{\prime },Q^{\prime }),
\end{equation*}%
where $(\beta :\delta ):(T:A)\rightarrow (T^{\prime }:A^{\prime })$ and $%
(\gamma :\sigma ):(S:B)\rightarrow (S^{\prime }:B^{\prime })$ are rng
morphisms, $\phi :P\rightarrow P^{\prime }$ is $(T,S)$-bilinear and $\psi
:Q\rightarrow Q^{\prime }$ is $(S,T)$-bilinear, such that%
\begin{equation*}
\beta (<p,q>_{T})=<\phi (p),\psi (q)>_{T^{\prime }}\text{ for all }p\in
P,q\in Q\text{ }.
\end{equation*}%
Notice that we consider $P^{\prime }$ as a $(T,S)$-bimodule and $Q^{\prime }$
as an $(S,T)$-bimodule with actions induced by the morphism of rngs $(\beta
:\delta )$ and $(\gamma :\sigma ).$ By $\mathbb{MSC}$ we denote the \emph{%
category of Morita semi-contexts} with morphisms defined as above, and by $%
\mathbb{UMSC}<\mathbb{MSC}$ the (non-full) subcategory of \emph{unital
Morita semi-contexts.}
\end{punto}

Morita semi-contexts are closely related to \emph{dual pairings} in the
sense of \cite{Abu-2005}:

\begin{punto}
Let $(T,S,P,Q,<,>_{T})\in \mathbb{MSC}$ and consider the canonical
isomorphisms of Abelian groups%
\begin{equation*}
\mathrm{Hom}_{(S,T)}(Q,\text{ }^{\ast }P)\overset{\xi }{\simeq }\text{ }%
\mathrm{Hom}_{(T,T)}(P\otimes _{S}Q,T)\overset{\zeta }{\simeq }\mathrm{Hom}%
_{(T,S)}(P,Q^{\ast }).
\end{equation*}%
This means that we have two dual $T$-pairings $\mathbf{P}_{l}:=(Q,$ $%
_{T}P)\in \mathcal{P}_{l}(T)$ and $\mathbf{Q}_{r}:=(P,Q_{T})\in \mathcal{P}%
_{r}(T),$ induced by the canonical $T$-linear morphisms%
\begin{equation*}
\kappa _{\mathbf{P}_{l}}:=\xi ^{-1}(<,>_{T}):Q\rightarrow \text{ }^{\ast }P%
\text{ and }\kappa _{\mathbf{Q}_{r}}:=\zeta (<,>_{T}):P\rightarrow Q^{\ast }.
\end{equation*}%
On the other hand, let $(S,T,Q,P,<,>_{S})\in \mathbb{MSC}$ and consider the
canonical isomorphisms of Abelian groups%
\begin{equation*}
\mathrm{Hom}_{(S,T)}(Q,P^{\ast })\overset{\xi ^{\prime }}{\simeq }\text{ }%
\mathrm{Hom}_{(S,S)}(Q\otimes _{T}P,S)\overset{\zeta ^{\prime }}{\simeq }%
\mathrm{Hom}_{(T,S)}(P,\text{ }^{\ast }Q).
\end{equation*}%
Then we have two dual $S$-pairings $\mathbf{P}_{r}:=(Q,P_{S})\in \mathcal{P}%
_{r}(S)$ and $\mathbf{Q}_{l}:=(P,$ $_{S}Q)\in \mathcal{P}_{l}(S),$ induced
by the canonical morphisms
\begin{equation*}
\kappa _{\mathbf{P}_{r}}:=\xi ^{\prime -1}(<,>_{S}):Q\rightarrow P^{\ast }%
\text{ and }\kappa _{\mathbf{Q}_{r}}:=\zeta ^{\prime }(<,>_{S}):P\rightarrow
\text{ }^{\ast }Q.
\end{equation*}
\end{punto}

\begin{punto}
By a \textbf{Morita datum }we mean a tuple%
\begin{equation}
\mathcal{M}=((T:A),(S:B),P,Q,<,>_{T},<,>_{S},I,J),  \label{Mori}
\end{equation}%
where the following are Morita semi-contexts.
\begin{equation}
\mathcal{M}_{T}:=((T:A),(S:B),P,Q,<,>_{T},I)\text{ and }\mathcal{M}%
_{S}:=((S:B),(T:A),Q,P,<,>_{S},J)  \label{MT-MS}
\end{equation}%
If, moreover, the bilinear morphisms $<,>_{T}:P\otimes _{S}Q\rightarrow T$
and $<-,>_{S}:Q\otimes _{T}P\rightarrow S$ are \emph{compatible, }in the
sense that%
\begin{equation}
<q,p>_{S}q^{\prime }=q<p,q^{\prime }>_{T}\text{ and }p<q,p^{\prime }>_{S}%
\text{ }=<p,q>_{T}p^{\prime }\text{ }\forall \text{ }p,p^{\prime }\in P,%
\text{ }q,q^{\prime }\in Q,  \label{M-compatible}
\end{equation}%
then we call $\mathcal{M}$ a \textbf{Morita context}. If $T,$ $S$ in a
Morita datum (context) $\mathcal{M}$ are unital, then we call $\mathcal{M}$
a \textbf{unital Morita datum (context).}
\end{punto}

\begin{punto}
Let $\mathcal{M}=((T:A),(S:B),P,Q,<,>_{T},<,>_{S})$ and $\mathcal{M}^{\prime
}=((T^{\prime }:A^{\prime }),(S^{\prime }:B^{\prime }),P^{\prime },Q^{\prime
},<,>_{T^{\prime }},<,>_{S^{\prime }})$ be Morita contexts. Extending \cite[%
Page 275]{Ami-1971}, we mean by a \textbf{morphism of Morita contexts} from $%
\mathcal{M}$ to $\mathcal{M}^{\prime }$ a four fold set of maps%
\begin{equation*}
((\beta :\delta ),(\gamma :\sigma ),\phi ,\psi
):((T:A),(S:B),P,Q)\rightarrow ((T^{\prime }:A^{\prime }),(S^{\prime
}:B^{\prime }),P^{\prime },Q^{\prime }),
\end{equation*}%
where $(\beta :\delta ):(T:A)\rightarrow (T^{\prime }:A^{\prime }),$ $%
(\gamma :\sigma ):(S:B)\rightarrow (S^{\prime }:B^{\prime })$ are rng
morphisms, $\phi :P\rightarrow P^{\prime }$ is $(T,S)$-bilinear and $\psi
:Q\rightarrow Q^{\prime }$ is $(S,T)$-bilinear, such that
\begin{equation*}
\beta (<p,q>_{T})=<\phi (p),\psi (q)>_{T^{\prime }}\text{ and }\gamma
(<q,p>_{S})=<\psi (q),\phi (p)>_{S^{\prime }}\text{ }\forall \text{ }p\in
P,q\in Q.\text{ }
\end{equation*}%
By $\mathbb{MC}$ we denote the \emph{category of Morita contexts} with
morphisms defined as above, and by $\mathbb{UMC}<\mathbb{MC}$ the (non-full)
subcategory of \emph{unital Morita contexts}.
\end{punto}

\begin{ex}
If $R$ is commutative, then any Morita semi-context $(R,R,P,Q,<,>_{R})$
yields a Morita context $(R,R,P,Q,<,>_{R},[,]_{R}),$ where $%
[,]_{R}:=Q\otimes _{R}P\simeq P\otimes _{R}Q\overset{<,>_{R}}{%
\longrightarrow }R.\blacksquare $
\end{ex}

\begin{punto}
We call a Morita semi-context $\mathbf{m}_{T}=(T,S,P,Q,<,>_{T})$ \textbf{%
semi-derived }(\textbf{derived}),\textbf{\ }iff $S:=\mathrm{End}(_{T}P)^{op}$
(and $Q=$ $^{\ast }P$). We call a Morita datum, or a Morita context, $%
\mathcal{M}=(T,S,P,Q,<,>_{T},<,>_{S})$ \textbf{semi-derived (derived)}, iff $%
S=\mathrm{End}(_{T}P)^{op},$ or $T=\mathrm{End}(P_{S})$ ($S=\mathrm{End}%
(_{T}P)^{op}$ and $Q=$ $^{\ast }P,$ or $T=\mathrm{End}(P_{S})$ and $%
Q=P^{\ast }$).
\end{punto}

\begin{remark}
Following \cite[1.2.]{Cae1998} (however, dropping the condition that the
bilinear map $<,>_{T}:P\otimes _{S}Q\rightarrow T$ is surjective), Morita
semi-contexts $(T,S,P,Q,<>_{T})$ in our sense were called \emph{dual pairs}
in \cite{Ver-2006}. However, we think the terminology we are using is more
informative and avoids confusion with other notions of dual pairings in the
literature (e.g. the ones studied by the first author in \cite{Abu-2005}).
The reason for this specific terminology (i.e. Morita semi-contexts) is that
every Morita context contains two Morita semi-contexts as clear from the
definition; and that any Morita semi-context can be \emph{extended} to a
(not necessarily unital) Morita context in a natural way as explained below.
\end{remark}

\subsection*{Elementary rngs}

\qquad In what follows we demonstrate how to build new Morita
(semi-)contexts from a given Morita semi-context. These constructions are
inspired by the notion of \emph{elementary rngs} in \cite[1.2.]{Cae1998}
(and \cite[Remark 3.8.]{Ver-2006}):

\begin{lemma}
\label{<>T}Let $\mathbf{m}_{T}:=((T:A),(S:B),P,Q,<,>_{T})\in \mathbb{MSC}.$

\begin{enumerate}
\item The $(T,T)$-bimodule $\mathbb{T}:=P\otimes _{S}Q$ has a structure of a
$T$-rng \emph{(}$A$-rng\emph{) }with multiplication%
\begin{equation*}
(p\otimes _{S}q)\cdot _{\mathbb{T}}(p^{\prime }\otimes _{S}q^{\prime
}):=<p,q>_{T}p^{\prime }\otimes _{S}\text{ }q^{\prime }\text{ }\forall \text{
}p,p^{\prime }\in P,\text{ }q,q^{\prime }\in Q,
\end{equation*}%
such that $<,>_{T}:\mathbb{T}\rightarrow T$ is a morphism of $A$-rngs, $P$
is a $(\mathbb{T},S)$-bimodule and $Q$ is an $(S,\mathbb{T})$-bimodule, where%
\begin{equation*}
(p\otimes _{S}q)\rightharpoonup \widetilde{p}:=<p,q>_{T}\widetilde{p}\text{
and }\widetilde{q}\leftharpoonup (p\otimes _{S}q):=\widetilde{q}<p,q>_{T}.
\end{equation*}%
Moreover, we have morphisms of $T$-rngs \emph{(}$A$-rngs\emph{)}%
\begin{equation*}
\begin{tabular}{llllllll}
$\psi $ & $:$ & $\mathbb{T}$ & $\rightarrow $ & $\mathrm{End}(P_{S}),$ & $%
p\otimes _{S}q$ & $\mapsto $ & $[\widetilde{p}\mapsto <p,q>_{T}\widetilde{p}%
];$ \\
$\phi $ & $:$ & $\mathbb{T}$ & $\rightarrow $ & $\mathrm{End}(_{S}Q)^{op},$
& $p\otimes _{S}q$ & $\mapsto $ & $[\widetilde{q}\mapsto \widetilde{q}%
<p,q>_{T}],$%
\end{tabular}%
\end{equation*}%
$((\mathbb{T}:A),(S:B),P,Q,id_{\mathbb{T}})\in \mathbb{MSC}$ and we have a
morphism of Morita semi-contexts%
\begin{equation*}
(<,>_{T},id_{S},,id_{P},id_{Q}):(\mathbb{T},S,P,Q,id_{\mathbb{T}%
})\rightarrow (T,S,P,Q,<,>_{T}).
\end{equation*}

\item The $(S,S)$-bimodule $\mathbf{S}:=Q\otimes _{T}P$ has a structure of
an $S$-rng \emph{(}$B$-rng\emph{) }with multiplication%
\begin{equation*}
(q\otimes _{T}p)\cdot _{\mathbf{S}}(q^{\prime }\otimes _{T}p^{\prime
}):=q<p,q^{\prime }>_{T}\otimes _{T}\text{ }p^{\prime }=q\otimes
_{T}<p,q^{\prime }>_{T}p^{\prime }\text{ }\forall \text{ }p,p^{\prime }\in P,%
\text{ }q,q^{\prime }\in Q,
\end{equation*}%
such that $<,>_{S}:\mathbf{S}\rightarrow S$ is a morphism of $B$-rngs, $P$
is a $(T,\mathbf{S})$-bimodule and $Q$ is an $(\mathbf{S},T)$-bimodule,
where
\begin{equation*}
\widetilde{p}\leftharpoonup (q\otimes _{T}p):=<\widetilde{p},q>_{T}p\text{
and }(q\otimes _{T}p)\rightharpoonup \widetilde{q}:=q<p,\widetilde{q}>_{T}.
\end{equation*}%
Moreover, we have morphisms of $S$-rngs \emph{(}$B$-rngs\emph{)}%
\begin{equation*}
\begin{tabular}{llllllll}
$\Psi $ & $:$ & $\mathbf{S}$ & $\rightarrow $ & $\mathrm{End}(_{T}P)^{op},$
& $q\otimes _{T}p$ & $\mapsto $ & $[\widetilde{p}\mapsto <\widetilde{p}%
,q>_{T}p],$ \\
$\Phi $ & $:$ & $\mathbf{S}$ & $\rightarrow $ & $\mathrm{End}(Q_{T}),$ & $%
q\otimes _{T}p$ & $\mapsto $ & $[\widetilde{q}\mapsto q<p,\widetilde{q}%
>_{T}],$%
\end{tabular}%
\end{equation*}%
and $\mathcal{M}:=((T:A),(\mathbf{S}:B),P,Q,<,>_{T},id_{\mathbf{S}})$ is a
Morita context.
\end{enumerate}
\end{lemma}

\begin{remarks}
\begin{enumerate}
\item Given $((S:B),(T:A),Q,P,<,>_{S})\in \mathbb{MSC},$ the $(S,S)$%
-bimodule $\mathbb{S}:=Q\otimes _{T}P$ becomes an $S$-rng with multiplication%
\begin{equation*}
(q\otimes _{T}p)\cdot _{\mathbb{S}}(q^{\prime }\otimes _{T}p^{\prime
}):=<q,p>_{S}q^{\prime }\otimes _{T}\text{ }p^{\prime }\text{ }\forall \text{
}p,p^{\prime }\in P,\text{ }q,q^{\prime }\in Q;
\end{equation*}%
and the $(T,T)$-bimodule $\mathbf{T}:=P\otimes _{S}Q$ becomes a $T$-rng with
multiplication%
\begin{equation*}
(p\otimes _{S}q)\cdot _{\mathbf{T}}(p^{\prime }\otimes _{S}q^{\prime
}):=p<q,p^{\prime }>_{S}\otimes _{S}\text{ }q^{\prime }=p\otimes _{S}\text{ }%
<q,p^{\prime }>_{S}q^{\prime }\text{ }\forall \text{ }p,p^{\prime }\in P,%
\text{ }q,q^{\prime }\in Q.
\end{equation*}%
Analogous results to those in Lemma \ref{<>T} can be obtained for the $S$%
-rng $\mathbb{S}$ and the $T$-rng $\mathbf{T}.$

\item Given a Morita semi-context $(T,S,P,Q,<,>_{T})$ several equivalent
conditions for the $T$-rng $\mathbf{T}:=P\otimes _{S}Q$ to be unital and the
modules $_{\mathbf{T}}P,$ $Q_{\mathbf{T}}$ to be \emph{firm} can be found in
\cite[Theorem 3.3.]{Ver-2006}. Analogous results can be formulated for the $%
S $-rng $Q\otimes _{T}P$ and the $S$-modules $P_{\mathbf{S}},$ $_{\mathbf{S}%
}Q$ corresponding to any $(S,T,Q,P,<,>_{S})\in \mathbb{MSC}.$
\end{enumerate}
\end{remarks}

\begin{proposition}
\begin{enumerate}
\item Let $\mathbf{m}_{T}=(T,S,P,Q,<,>_{T})\in \mathbb{UMSC}$ and assume the
$A$-rng $\mathbb{T}:=P\otimes _{S}Q$ to be unital. If $<,>_{T}:\mathbb{T}%
\rightarrow T$ respects unities \emph{(}and $\mathbf{m}_{T}$ is injective%
\emph{),} then $<,>_{T}$ is surjective \emph{(}$\mathbb{T}\overset{<,>_{T}}{%
\simeq }T$ as $A$-rings\emph{)}.

\item Let $\mathbf{m}_{S}=(S,T,Q,P,<,>_{S})\in \mathbb{UMSC}$ and assume the
$B$-rng $\mathbb{S}:=Q\otimes _{S}P$ to be unital. If $<,>_{S}:\mathbb{S}%
\rightarrow S$ respects unities \emph{(}and $\mathbf{m}_{S}$ is injective%
\emph{)}, then $<,>_{S}$ is surjective \emph{(}$\mathbb{S}\overset{<,>_{S}}{%
\simeq }S$ as $B$-rings\emph{)}.

\item Let $\mathcal{M}=(T,S,P,Q,<,>_{T},<,>_{S})\in \mathbb{UMC}$ and assume
the rngs $\mathbb{T}:=P\otimes _{S}Q,$ $T,$ $\mathbb{S}:=Q\otimes _{S}P$ to
be unital. If $<,>_{T}:$ $P\otimes _{S}Q\rightarrow T$ and $<,>_{S}:\mathbb{S%
}\rightarrow S$ respect unities, then $\mathbb{T}\overset{<,>_{T}}{\simeq }T$
as $A$-ring, $\mathbb{S}\overset{<,>_{S}}{\simeq }S$ as $B$-rings and we
have equivalences of categories $_{\mathbb{T}}\mathbb{M}\approx $ $_{\mathbb{%
S}}\mathbb{M}$ \emph{(}and $\mathbb{M}_{\mathbb{T}}\approx \mathbb{M}_{%
\mathbb{S}}$\emph{)}.
\end{enumerate}
\end{proposition}

\begin{Beweis}
Assume $\mathbb{T}$ is unital with $1_{\mathbb{T}}=\dsum_{i=1}^{n}p_{i}%
\otimes _{S}q_{i}.$ If $<,>_{T}$ respects unities, then we have $%
\dsum_{i=1}^{n}<p_{i},q_{i}>_{T}=1_{T},$ and so for any $t\in T$ we get $%
t=t1_{T}=\dsum_{i=1}^{n}t<p_{i},q_{i}>_{T}=\dsum_{i=1}^{n}<tp_{i},q_{i}>_{T}%
\in \mathrm{Im}(<,>_{T}).$ One can prove \textquotedblleft
2\textquotedblright\ analogously. As for \textquotedblleft
3\textquotedblright , it is well known that a unital Morita context with
surjective connecting bimodule morphisms is strict (e.g. \cite[12.7.]%
{Fai-1981}), hence $\mathbb{T}\overset{<,>_{T}}{\simeq }T,$ $\mathbb{S}%
\overset{<,>_{S}}{\simeq }S.$ The equivalences of categories $_{\mathbb{T}}%
\mathbb{M}\simeq $ $_{T}\mathbb{M}\approx $ $_{S}\mathbb{M}\simeq $ $_{%
\mathbb{S}}\mathbb{M}$ (and $\mathbb{M}_{\mathbb{T}}\simeq \mathbb{M}%
_{T}\approx \mathbb{M}_{S}\simeq \mathbb{M}_{\mathbb{S}}$)\emph{\ }follow
then by classical Morita Theory (e.g. \cite[Chapter 12]{Fai-1981}).$%
\blacksquare $
\end{Beweis}

\begin{definition}
Let $T$ be an $A$-rng, $V_{T}$ a right $T$-module and consider for every
left $T$-module $_{T}L$ the annihilator
\begin{equation*}
\mathrm{ann}_{L}^{\otimes }(V_{T}):=\{l\in L\mid V\otimes _{T}l=0\}.
\end{equation*}%
Following \cite[Exercises 19]{AF-1974}, we say $V_{T}$ is $L$\textbf{%
-faithful}, iff $\mathrm{ann}_{L}^{\otimes }(V_{T})=0;$ and to be \textbf{%
completely faithful}, iff $V_{T}$ is $L$-faithful for every left $T$-module $%
_{S}L.$ Similarly, we can define completely faithful left $T$-modules.
\end{definition}

\qquad Under suitable conditions, the following result characterizes the
Morita data, which are Morita contexts:

\begin{proposition}
\label{T=T}Let $\mathcal{M}=(T,S,P,Q,<,>_{T},<,>_{S})$ be a Morita datum.

\begin{enumerate}
\item If $\mathcal{M}\in \mathbb{MC}$, then $\mathbf{S}\overset{id}{\simeq }%
\mathbb{S}$ and $\mathbf{T}\overset{id}{\simeq }\mathbb{T}$ as rngs.

\item Assume $_{T}P$ is $Q$-faithful and $Q_{T}$ is $P$-faithful. Then $%
\mathcal{M}\in \mathbb{MC}$ if and only if $\mathbf{S}\overset{id}{\simeq }%
\mathbb{S}$ and $\mathbf{T}\overset{id}{\simeq }\mathbb{T}$ as rngs.
\end{enumerate}
\end{proposition}

\begin{Beweis}
\begin{enumerate}
\item Obvious.

\item Assume $\mathbf{S}\overset{id}{\simeq }\mathbb{S}$ and $\mathbf{T}%
\overset{id}{\simeq }\mathbb{T}$ as rngs. If $p\in P$ and $q,q^{\prime }\in
Q $ are arbitrary, then we have for any $\widetilde{p}\in P:$%
\begin{equation*}
<q,p>_{S}q^{\prime }\otimes _{T}\text{ }\widetilde{p}=(q\otimes _{T}p)\cdot
_{\mathbb{S}}(q^{\prime }\otimes _{T}\widetilde{p})=(q\otimes _{T}p)\cdot _{%
\mathbf{S}}(q^{\prime }\otimes _{T}\widetilde{p})=q<p,q^{\prime
}>_{T}\otimes _{T}\text{ }\widetilde{p},
\end{equation*}%
hence $<q,p>_{S}q^{\prime }-q<p,q^{\prime }>_{T}\in \mathrm{ann}_{Q}(P)=0$
(since $_{T}P$ is $Q$-faithful), i.e. $<q,p>_{S}q^{\prime }=q<p,q^{\prime
}>_{T}$ for all $p\in P$ and $q,q^{\prime }\in Q.$ Assuming $Q_{T}$ is $P$%
-faithful, one can prove analogously that $<p,q>_{T}p^{\prime
}=p<q,p^{\prime }>_{S}$ for all $p,p^{\prime }\in P$ and $q\in Q.$
Consequently, $\mathcal{M}$ is a Morita context.$\blacksquare $
\end{enumerate}
\end{Beweis}

\section{Injective Morita (Semi-)Contexts}

\begin{definition}
We call a Morita semi-context $\mathbf{m}_{T}=(T,S,P,Q,<,>_{T},I):$

\textbf{injective} (resp. \textbf{semi-strict}, \textbf{strict}), iff $%
<,>_{T}:P\otimes _{S}Q\rightarrow T$ is injective (resp. surjective,
bijective);

\textbf{non-degenerate}, iff $Q\hookrightarrow $ $^{\ast }P$ and $%
P\hookrightarrow Q^{\ast }$ canonically;

\textbf{Morita }$\alpha $\textbf{-semi-context}, iff $\mathbf{P}_{l}:=(Q,$ $%
_{T}P)\in \mathcal{P}_{l}^{\alpha }(T)$ and $\mathbf{Q}_{r}:=(P,Q_{T})\in
\mathcal{P}_{r}^{\alpha }(T).$
\end{definition}

\begin{notation}
By $\mathbb{MSC}^{\alpha }\leq \mathbb{MSC}$ ($\mathbb{UMSC}^{\alpha }\leq
\mathbb{UMSC}$) we denote the full subcategory of (unital) Morita
semi-contexts satisfying the $\alpha $-condition. Moreover, we denote by $%
\mathbb{IMSC}\leq \mathbb{MSC}$ ($\mathbb{IUMSC}\leq \mathbb{UMSC}$) the
full subcategory of injective (unital) Morita semi-contexts.
\end{notation}

\begin{definition}
We say a Morita datum (context) $\mathcal{M}=(T,S,P,Q,<,>_{T},<,>_{S},I,J):$

is \textbf{injective} (resp. \textbf{semi-strict}, \textbf{strict}), iff $%
<,>_{T}:P\otimes _{S}Q\rightarrow T$ and $<,>_{S}:Q\otimes _{T}P\rightarrow
S $ are injective (resp. surjective, bijective);

is \textbf{non-degenerate}, iff $Q\hookrightarrow $ $^{\ast }P,$ $%
P\hookrightarrow Q^{\ast },$ $Q\hookrightarrow P^{\ast }$ and $%
P\hookrightarrow $ $^{\ast }Q$ canonically;

\textbf{satisfies the left }$\alpha $\textbf{-condition}, iff $\mathbf{P}%
_{l}:=(Q,$ $_{T}P)\in \mathcal{P}_{l}^{\alpha }(T)$ and $\mathbf{Q}_{l}:=(P,$
$_{S}Q)\in \mathcal{P}_{l}^{\alpha }(S);$

\textbf{satisfies the right }$\alpha $\textbf{-condition,} iff $\mathbf{Q}%
_{r}:=(P,Q_{T})\in \mathcal{P}_{r}^{\alpha }(T)$ and $\mathbf{P}%
_{r}:=(Q,P_{S})\in \mathcal{P}_{r}^{\alpha }(S);$

\textbf{satisfies the }$\alpha $\textbf{-condition,} or $\mathcal{M}$ is a
\textbf{Morita }$\alpha $\textbf{-datum (Morita }$\alpha $-\textbf{context}%
), iff $\mathcal{M}$ satisfies both the left and the right $\alpha $%
-conditions.
\end{definition}

\begin{notation}
By $\mathbb{MC}_{l}^{\alpha }<\mathbb{MC}$ ($\mathbb{UMC}_{l}^{\alpha }<%
\mathbb{UMC}$) we denote the full subcategory of Morita contexts satisfying
the left $\alpha $-condition, and by $\mathbb{MC}_{r}^{\alpha }<\mathbb{MC}$
($\mathbb{UMC}_{r}^{\alpha }<\mathbb{UMC}$) the full subcategory of (unital)
Morita contexts satisfying the right $\alpha $-condition. Moreover, we set $%
\mathbb{MC}^{\alpha }:=\mathbb{MC}_{l}^{\alpha }\cap \mathbb{MC}_{r}^{\alpha
}$ and $\mathbb{UMC}^{\alpha }:=\mathbb{UMC}_{l}^{\alpha }\cap \mathbb{UMC}%
_{r}^{\alpha }.$
\end{notation}

\begin{lemma}
\label{alpha->inj}Let $\mathcal{M}=(T,S,P,Q,<,>_{T},<,>_{S},I,J)\in \mathbb{%
MC}.$ Consider the Morita semi-context $\mathcal{M}_{S}:=(S,T,Q,P,<,>_{S}),$
the dual pairings $\mathbf{P}_{l}:=(Q,$ $_{T}P)\in \mathcal{P}_{l}(T),$ $%
\mathbf{Q}_{r}:=(P,Q_{T})\in \mathcal{P}_{r}(T)$ and the canonical morphisms
of rings
\begin{equation*}
\rho _{P}:S\rightarrow \mathrm{End}(_{T}P)^{op}\text{ and }\lambda
_{Q}:S\rightarrow \mathrm{End}(Q_{T}).
\end{equation*}

\begin{enumerate}
\item If $\mathbf{Q}_{r}$ is injective \emph{(}semi-strict\emph{)}, then $%
\mathcal{M}_{S}$ is injective \emph{(}$\rho _{P}:S\rightarrow \mathrm{End}%
(_{T}P)^{op}$ is a surjective morphism of $B$-rngs\emph{)}.

\item Assume $P_{S}$ is faithful and let $\mathbf{Q}_{r}$ be semi-strict.
Then $S\simeq \mathrm{End}(_{T}P)^{op}$ \emph{(}an isomorphism of unital $B$%
-rings\emph{) }and $\mathcal{M}_{S}$ is strict.

\item If $\mathbf{P}_{l}$ is injective \emph{(}semi-strict\emph{)}, then $%
\mathcal{M}_{S}$ is injective \emph{(}$\lambda _{Q}:S\rightarrow \mathrm{End}%
(Q_{T})$ is a surjective morphism of $B$-rngs\emph{)}.

\item Assume $_{S}Q$ is faithful and let $\mathbf{P}_{l}$ is semi-strict.
Then $S\simeq \mathrm{End}(Q_{T})$ \emph{(}an isomorphism of unital $B$-rings%
$\emph{)}$ and $\mathcal{M}_{S}$ is strict.
\end{enumerate}
\end{lemma}

\begin{Beweis}
We prove only \textquotedblleft 1\textquotedblright\ and \textquotedblleft
2\textquotedblright , as \textquotedblleft 3\textquotedblright\ and
\textquotedblleft 4\textquotedblright\ can be proved analogously.

Consider the following butterfly diagram with canonical morphisms
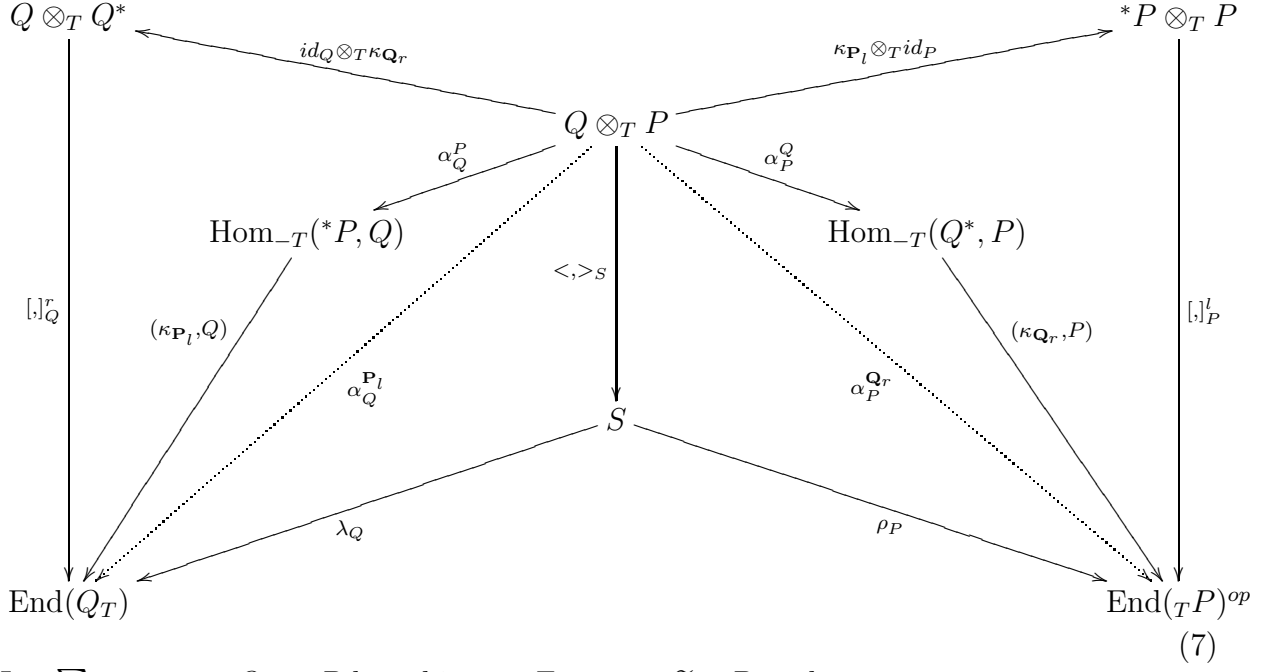
\begin{equation}
\xymatrix{ Q \otimes_T Q^{*} \ar[dddddd]_{[,]_Q ^{r}} & & & & & & {}^{*}P
\otimes_T P \ar[dddddd]^{[,]_P ^{l}} \\ & & & Q \otimes_{T} P
\ar[ulll]_{id_Q \otimes_T \kappa_{\mathbf{Q}_r}}
\ar[urrr]^{\kappa_{\mathbf{P}_l} \otimes_T id_P} \ar[ddd]_{<,>_S}
\ar[rrd]^{\alpha_P^Q} \ar[lld]_{\alpha_Q^P}
\ar@{.>}[rrrddddd]_{\alpha_P^{\mathbf{Q}_r}}
\ar@{.>}[dddddlll]^{\alpha_Q^{\mathbf{P}_l}}& & & \\ & {\rm
Hom}_{-T}(^{*}P,Q) \ar[ddddl]_(0.3){(\kappa_{\mathbf{P}_l},Q)} & & & & {\rm
Hom}_{-T}(Q^{*},P) \ar[ddddr]^(0.3){(\kappa_{\mathbf{Q}_r},P)} & \\ & & & &
& & \\ & & & S \ar[ddlll]^{\lambda_Q} \ar[ddrrr]_{\rho_P} & & & \\ & & & & &
& \\ {\rm End}(Q_{T}) & & & & & & {\rm End}(_{T} P)^{op}}  \label{S1}
\end{equation}%
Let $\sum q_{i}\otimes _{T}p_{i}\in Q\otimes _{T}P$ be arbitrary. For every $%
\widetilde{p}\in P$ we have
\begin{equation*}
\begin{tabular}{lll}
$\lbrack (\kappa _{\mathbf{Q}_{r}},P)\circ \alpha _{P}^{Q})(\sum
q_{i}\otimes _{T}p_{i})](\widetilde{p})$ & $=$ & $\sum <\widetilde{p}%
,q_{i}>_{T}p_{i}$ \\
& $=$ & $\sum \widetilde{p}<q_{i},p_{i}>_{S}$ \\
& $=$ & $\rho _{P}(\sum <q_{i},p_{i}>_{S})(\widetilde{p})$ \\
& $=$ & $(\rho _{P}\circ <,>_{S})(\sum q_{i}\otimes _{T}p_{i})(\widetilde{p}%
),$%
\end{tabular}%
\end{equation*}%
i.e. $\alpha _{P}^{\mathbf{Q}_{r}}:=(\kappa _{\mathbf{Q}_{r}},P)\circ \alpha
_{P}^{Q}=\rho _{P}\circ <,>_{S};$ and
\begin{equation*}
\begin{tabular}{lll}
$\lbrack ,]_{P}^{l}\circ (\kappa _{\mathbf{P}_{l}}\otimes _{T}id_{P}))(\sum
q_{i}\otimes _{T}p_{i})](\widetilde{p})$ & $=$ & $\dsum \kappa _{\mathbf{P}%
_{l}}(q_{i})(\widetilde{p})p_{i}$ \\
& $=$ & $\dsum <\widetilde{p},q_{i}>_{T}p_{i}$ \\
& $=$ & $\dsum \widetilde{p}<q_{i},p_{i}>_{S}$ \\
& $=$ & $\rho _{P}(\dsum <q_{i},p_{i}>_{S})(\widetilde{p})$ \\
& $=$ & $[(\rho _{P}\circ <,>_{S})(\sum q_{i}\otimes _{T}p_{i})](\widetilde{p%
}),$%
\end{tabular}%
\end{equation*}%
i.e. $[,]_{P}^{l}\circ (\kappa _{\mathbf{P}_{l}}\otimes _{T}id_{P})=\rho
_{P}\circ <,>_{S}.$ On the other hand, for every $\widetilde{q}\in Q$ we
have
\begin{equation*}
\begin{tabular}{lll}
$((\kappa _{\mathbf{P}_{l}},Q)\circ \alpha _{Q}^{\mathbf{P}_{l}})(\sum
q_{i}\otimes _{T}p_{i})(\widetilde{q})$ & $=$ & $\sum q_{i}<p_{i},\widetilde{%
q}>_{T}$ \\
& $=$ & $(\sum <q_{i},p_{i}>_{S})\widetilde{q}$ \\
& $=$ & $\lambda _{Q}(\sum <q_{i},p_{i}>_{S})(\widetilde{q})$ \\
& $=$ & $(\lambda _{Q}\circ <,>_{S})(\sum q_{i}\otimes _{T}p_{i}),$%
\end{tabular}%
\end{equation*}%
i.e. $\alpha _{Q}^{\mathbf{P}_{l}}:=(\kappa _{\mathbf{P}_{l}},Q)\circ \alpha
_{Q}^{\mathbf{P}_{l}}=\lambda _{Q}\circ <,>_{S}$ and
\begin{equation*}
\begin{tabular}{lll}
$([,]_{Q}^{r}\circ (id_{Q}\otimes _{T}\kappa _{\mathbf{Q}_{r}}))(\sum
q_{i}\otimes _{T}p_{i})](\widetilde{q})$ & $=$ & $\dsum q_{i}\kappa _{%
\mathbf{Q}_{r}}(p_{i})(\widetilde{q})$ \\
& $=$ & $\dsum q_{i}<p_{i},\widetilde{q}>_{T}$ \\
& $=$ & $\dsum <q_{i},p_{i}>_{S}\widetilde{q}$ \\
& $=$ & $\lambda _{Q}(\dsum <q_{i},p_{i}>_{S})(\widetilde{q})$ \\
& $=$ & $[(\lambda _{Q}\circ <,>_{S})(\sum q_{i}\otimes _{T}p_{i})](%
\widetilde{q}),$%
\end{tabular}%
\end{equation*}%
i.e. $[,]_{Q}^{r}\circ (id_{Q}\otimes _{T}\kappa _{\mathbf{Q}_{r}})=\lambda
_{Q}\circ <,>_{S}.$ Hence Diagram (\ref{S1}) is commutative.

(1) Follows directly from the assumptions and the equality $\alpha _{P}^{%
\mathbf{Q}_{r}}=\rho _{P}\circ <,>_{S}.$

(2) Let $P_{S}$ be faithful, so that the canonical left $S$-linear map $\rho
_{P}:S\rightarrow \mathrm{End}(_{T}P)^{op}$ is injective. Assume now that $%
\mathbf{Q}_{r}$ is semi-strict. Then $\rho _{P}$ is surjective by
\textquotedblleft 1\textquotedblright\ , whence bijective. Since rings of
endomorphisms are unital, we conclude that $S\simeq \mathrm{End}(_{T}P)^{op}$
is a unital $B$-ring as well (with unity $\rho _{P}^{-1}(id_{P})$).
Moreover, the surjectivity of $\alpha _{P}^{\mathbf{Q}_{r}}=\rho _{P}\circ
<,>_{S}$ implies that $<,>_{S}$ is surjective (since $\rho _{P}$ is
injective), say $1_{S}=\dsum_{j}<\widetilde{q}_{j},\widetilde{p}_{j}>_{S}$
for some $\{(\widetilde{q}_{j},\widetilde{p}_{j})\}_{J}\subseteq Q\times P.$
For any $\dsum_{i}q_{i}\otimes _{T}p_{i}\in \mathrm{Ker}(<,>_{S}),$ we have
then%
\begin{equation*}
\begin{tabular}{lllll}
$\dsum_{i}q_{i}\otimes _{T}p_{i}$ & $=$ & $(\dsum_{i}q_{i}\otimes
_{T}p_{i})\cdot 1_{S}$ & $=$ & $\dsum_{i}(q_{i}\otimes _{T}p_{i})\cdot
(\dsum_{j}<\widetilde{q}_{j},\widetilde{p}_{j}>_{S})$ \\
& $=$ & $\dsum_{i,j}q_{i}\otimes _{T}p_{i}<\widetilde{q}_{j},\widetilde{p}%
_{j}>_{S}$ & $=$ & $\dsum_{i,j}q_{i}\otimes _{T}<p_{i},\widetilde{q}_{j}>_{T}%
\widetilde{p}_{j}$ \\
& $=$ & $\dsum_{i,j}q_{i}<p_{i},\widetilde{q}_{j}>_{T}\otimes _{T}\widetilde{%
p}_{j}$ & $=$ & $\dsum_{i,j}<q_{i},p_{i}>_{S}\widetilde{q}_{j}\otimes _{T}%
\widetilde{p}_{j}$ \\
& $=$ & $\dsum_{j}(\dsum_{i}<q_{i},p_{i}>_{S})\widetilde{q}_{j}\otimes _{T}%
\widetilde{p}_{j}$ & $=$ & $0,$%
\end{tabular}%
\end{equation*}%
i.e. $<,>_{S}$ is injective, whence an isomorphism.$\blacksquare $
\end{Beweis}

\qquad

\qquad The following result shows that Morita $\alpha $-contexts are
injective:

\begin{corollary}
\label{LR-alpha-M}$\mathbb{MC}_{l}^{\alpha }\cup \mathbb{MC}_{r}^{\alpha
}\leq \mathbb{IMC}.$
\end{corollary}

\begin{ex}
\label{QF}Let $\mathbf{m}_{T}=(T,S,P,Q,<,>_{T})$ be a non-degenerate Morita
semi-context. If $T$ is a QF ring and the $T$-modules $_{T}P,$ $Q_{T}$ are
projective, then by Proposition \ref{rp-rp} \textquotedblleft
7\textquotedblright\ $\mathbf{P}_{l}:=(Q,$ $_{T}P)\in \mathcal{P}%
_{l}^{\alpha }(T)$ and $\mathbf{Q}_{r}:=(P,Q_{T})\in \mathcal{P}_{r}^{\alpha
}(T)$ (i.e. $\mathbf{m}_{T}$ is a Morita $\alpha $-semi-context, whence
injective). On the other hand, let $\mathcal{M}=(T,S,P,Q,<,>_{T},<,>_{S})$
be a non-degenerate Morita datum. If $T,S$ are QF rings and the modules $%
_{T}P,$ $Q_{T},$ $P_{S},$ $_{S}Q$ are projective, then $\mathcal{M}$ is an
Morita $\alpha $-datum (whence injective).$\blacksquare $
\end{ex}

Every semi-strict \emph{unital} Morita context is injective (whence strict,
e.g. \cite[12.7.]{Fai-1981}). The following example, which is a modification
of \cite[Example 18.30]{Lam-1999}), shows that the converse is not
necessarily true:

\begin{ex}
\label{iNs}Let $T=\mathrm{M}_{2}(\mathbb{Z}_{2})$ be the ring of $2\times 2$
matrices with entries in $\mathbb{Z}_{2}.$ Notice that $e=\left[
\begin{array}{cc}
1 & 0 \\
0 & 0%
\end{array}%
\right] \in T$ is an idempotent, and that $eTe\simeq \mathbb{Z}_{2}$ as
rings. Set%
\begin{equation*}
P:=Te=\{\left[
\begin{array}{cc}
a^{\prime } & 0 \\
c^{\prime } & 0%
\end{array}%
\right] \mid a^{\prime },c^{\prime }\in \mathbb{Z}_{2}\}\text{ and }Q:=eT=\{%
\left[
\begin{array}{cc}
a & b \\
0 & 0%
\end{array}%
\right] \mid a,b\in \mathbb{Z}_{2}\}.
\end{equation*}%
Then $P=Te$ is a $(T,eTe)$-bimodule and $Q=eT$ is an $(eTe,T)$-bimodule.
Moreover, we have a Morita context%
\begin{equation*}
\mathcal{M}_{e}=(T,eTe,Te,,eT,<,>_{T},<.>_{eTe}),
\end{equation*}%
where the connecting bilinear maps are%
\begin{equation*}
\begin{tabular}{lllll}
$<,>_{T}$ & $:$ & $Te\otimes _{eTe}eT$ & $\rightarrow $ & $T,$ \\
&  &  &  &  \\
&  & $\left[
\begin{array}{cc}
a^{\prime } & 0 \\
c^{\prime } & 0%
\end{array}%
\right] \otimes _{eTe}\left[
\begin{array}{cc}
a & b \\
0 & 0%
\end{array}%
\right] $ & $\mapsto $ & $\left[
\begin{array}{cc}
a^{\prime }a & a^{\prime }b \\
c^{\prime }a & c^{\prime }b%
\end{array}%
\right] $ \\
&  &  &  &  \\
$<,>_{eTe}$ & $:$ & $eT\otimes _{T}Te$ & $\rightarrow $ & $eTe$ \\
&  &  &  &  \\
&  & $\left[
\begin{array}{cc}
a & b \\
0 & 0%
\end{array}%
\right] \otimes _{T}\left[
\begin{array}{cc}
a^{\prime } & 0 \\
c^{\prime } & 0%
\end{array}%
\right] $ & $\mapsto $ & $\left[
\begin{array}{cc}
aa^{\prime }+bc^{\prime } & 0 \\
0 & 0%
\end{array}%
\right] .$%
\end{tabular}%
\end{equation*}%
Straightforward computations show that $<,>_{T}$ is injective but not
surjective (as $\left[
\begin{array}{cc}
1 & 1 \\
1 & 0%
\end{array}%
\right] \notin \mathrm{Im}(<,>_{T})$) and that $<,>_{eTe}$ is in fact an
isomorphism. This means that $\mathcal{M}_{e}$ is an injective Morita
context that is not semi-strict (whence not strict).$\blacksquare $
\end{ex}

\begin{definition}
\label{SF}Let $T$ be a rng and $I\vartriangleleft T$ an ideal. For every
left $T$-module $_{T}V$ consider the canonical $T$-linear map%
\begin{equation*}
\zeta _{I,V}:V\rightarrow \mathrm{Hom}_{T}(I,V),\text{ }v\mapsto \lbrack
t\mapsto tv].
\end{equation*}%
We say $_{T}I$ is \textbf{strongly }$V$\textbf{-faithful}, iff $\mathrm{ann}%
_{V}(I):=\mathrm{Ker}(\zeta _{I,V}):=0.$ Moreover, we say $I$ is \textbf{%
strongly faithful}, if $_{T}I$ is $V$-faithful for every left $T$-module $%
_{T}V.$ Strong faithfulness of $I$ w.r.t. right $T$-modules can be defined
analogously.
\end{definition}

\begin{remark}
Let $T$ be a rng, $I\vartriangleleft T$ an ideal and $_{T}U$ a left ideal.
It's clear that $\mathrm{ann}_{U}^{\otimes }(I_{T})\subseteq \mathrm{ann}%
_{U}(I):=\mathrm{Ker}(\zeta _{I,U}).$ Hence, if $_{T}I$ is \emph{strongly }$%
U $\emph{-faithful}, then $I_{T}$ is $U$\emph{-faithful} (which justifies
our terminology). In particular, if $_{T}I$ is \emph{strongly faithful},
then $I_{T}$ is \emph{completely faithful}.
\end{remark}

\qquad Morita $\alpha $-contexts are injective by Corollary \ref{LR-alpha-M}%
. The following result gives a partial converse:

\begin{lemma}
\label{part-conv}Let $\mathcal{M}=(T,S,P,Q,<,>_{T},<,>_{S},I,J)\in \mathbb{MC%
}$ and assume the Morita semi-context $\mathcal{M}_{S}:=(S,T,Q,P,<,>_{S},J)$
is injective.

\begin{enumerate}
\item If $_{S}J$ is strongly faithful, then $\mathbf{Q}_{r}:=(P,Q_{T})\in
\mathcal{P}_{r}^{\alpha }(T).$

\item If $J_{S}$ is strongly faithful, then $\mathbf{P}_{l}:=(Q,$ $_{T}P)\in
\mathcal{P}_{l}^{\alpha }(T).$
\end{enumerate}
\end{lemma}

\begin{Beweis}
We prove only \textquotedblleft 1\textquotedblright , since
\textquotedblleft 2\textquotedblright\ can be proved similarly. Assume $%
\mathcal{M}_{S}$ is injective and consider for every left $T$-module $U$ the
following diagram%
\begin{equation}
\xymatrix{ Q \otimes_{T} U \ar[rr]^{\alpha _{U}^{\mathbf{Q}_{r}}}
\ar[dr]_{\zeta_{J,Q \otimes_T U}} & & {\rm Hom}_{T-}(P,U)
\ar[dl]^{\psi_{Q,U}} \\ & {\rm Hom_{S-}}(J,Q \otimes_T U) & }  \label{J1}
\end{equation}%
where for all $f\in \mathrm{Hom}_{T-}(P,U)$ and $\dsum <q_{j},p_{j}>_{S}\in
J $ we define%
\begin{equation*}
\psi _{Q,U}(f)(\dsum <q_{j},p_{j}>_{S}):=\dsum q_{j}\otimes _{T}f(p_{j}).
\end{equation*}%
Then we have for every $\dsum \widetilde{q}_{i}\otimes _{T}\widetilde{u}%
_{i}\in Q\otimes _{T}U$ and $s=\dsum_{j}<q_{j},p_{j}>_{S}\in J:$%
\begin{equation*}
\begin{tabular}{lll}
$(\psi _{Q,U}\circ \alpha _{U}^{\mathbf{Q}_{r}})(\sum_{i}\widetilde{q}%
_{i}\otimes _{T}\widetilde{u}_{i})(s)$ & $=$ & $\dsum_{j}q_{j}\otimes
_{T}[\alpha _{U}^{\mathbf{Q}_{r}}(\dsum_{i}\widetilde{q}_{i}\otimes _{T}%
\widetilde{u}_{i})](p_{j})$ \\
& $=$ & $\dsum_{j}q_{j}\otimes _{T}\dsum_{i}<p_{j},\widetilde{q}_{i}>_{T}%
\widetilde{u}_{i}]$ \\
& $=$ & $\dsum_{i,j}q_{j}\otimes _{T}<p_{j},\widetilde{q}_{i}>_{T}\widetilde{%
u}_{i}$ \\
& $=$ & $\dsum_{i,j}q_{j}<p_{j},\widetilde{q}_{i}>_{T}\otimes _{T}\widetilde{%
u}_{i}$ \\
& $=$ & $\dsum_{i,j}<q_{j},p_{j}>_{S}\widetilde{q}_{i}\otimes _{T}\widetilde{%
u}_{i}$ \\
& $=$ & $\zeta _{J,Q\otimes _{T}U}(\sum_{i}\widetilde{q}_{i}\otimes _{T}%
\widetilde{u}_{i})(s),$%
\end{tabular}%
\end{equation*}%
i.e. diagram (\ref{J1}) is commutative. If $_{S}J$ is strongly faithful,
then $\mathrm{Ker}(\zeta _{J,Q\otimes _{T}U})=\mathrm{ann}_{Q\otimes
_{T}U}(J)=0,$ hence $\zeta _{J,Q\otimes _{T}U}$ is injective and it follows
then that $\alpha _{U}^{\mathbf{Q}_{r}}$ is injective.$\blacksquare $
\end{Beweis}

\begin{proposition}
Let $\mathcal{M}=(T,S,P,Q,<,>_{T},<,>_{S},I,J)\in \mathbb{IMC}.$ If $_{T}I,$
$I_{T},$ $_{S}J$ and $J_{S}$ are strongly faithful, then $\mathcal{M}\in
\mathbb{MC}^{\alpha }.$
\end{proposition}

\section{Equivalences of Categories}

\qquad In this section we give some applications of \emph{injective Morita }(%
\emph{semi-})\emph{contexts} and\emph{\ injective Morita data }to
equivalences between suitable subcategories of modules arising in the Kato-M%
\"{u}ller-Ohtake localization-colocalization theory (as developed in (e.g.
\cite{Kat-1978}, \cite{KO-1979}, \cite{Mul-1974}). All rings, hence all
Morita (semi-)contexts and data, in this section are unital.

\subsection*{Static and Adstatic Modules}

\begin{punto}
(\cite{C-IG-TW2003})\ Let $\mathcal{A}$ and $\mathcal{B}$ be two complete
cocomplete Abelian categories, $\mathbf{R}:\mathcal{A}\rightarrow \mathcal{B}
$ an additive covariant functor with left adjoint $\mathbf{L}:\mathcal{B}%
\rightarrow \mathcal{A}$ and let
\begin{equation*}
\omega :\mathbf{LR}\rightarrow 1_{\mathcal{A}}\text{ and }\eta :1_{\mathcal{B%
}}\rightarrow \mathbf{RL}
\end{equation*}%
be the induced natural transformations (called the \emph{counit} and the
\emph{unit} of the adjunction, respectively). Related to the adjoint pair $(%
\mathbf{L},\mathbf{R})$ are two \emph{full} subcategories of $\mathcal{A}$
and $\mathcal{B}:$%
\begin{equation*}
\mathrm{Stat}(\mathbf{R}):=\{X\in \mathcal{A}\mid \mathbf{LR}(X)\overset{%
\omega _{X}}{\simeq }X\}\text{ and }\mathrm{Adstat}(\mathbf{R}):=\{Y\in
\mathcal{B}\mid Y\overset{\eta _{Y}}{\simeq }\mathbf{RL}(Y)\},
\end{equation*}%
whose members are called $\mathbf{R}$-\textbf{static objects} and $\mathbf{R}
$-\textbf{adstatic objects,} respectively. It is evident (from definition)
that we have equivalence of categories $\mathrm{Stat}(\mathbf{R})\approx
\mathrm{Adstat}(\mathbf{R}).$
\end{punto}

\qquad A typical situation, in which static and adstatic objects arise
naturally is the following:

\begin{punto}
\label{ref=ref}Let $T,S$ be rings, $_{T}U_{S}$ a $(T,S)$-bimodule and
consider the covariant functors%
\begin{equation*}
\mathbf{H}_{U}^{l}:=\mathrm{Hom}_{T}(U,-):\text{ }_{T}\mathbb{M}\rightarrow
\text{ }_{S}\mathbb{M}\text{ and }\mathbf{T}_{U}^{l}:=U\otimes _{S}-:\text{ }%
_{S}\mathbb{M}\rightarrow \text{ }_{T}\mathbb{M}.
\end{equation*}%
It is well-known that $(\mathbf{T}_{U}^{l},\mathbf{H}_{U}^{l})$ is an
adjoint pair of covariant functors via the \emph{natural isomorphisms}%
\begin{equation*}
\mathrm{Hom}_{T}(U\otimes _{S}M,N)\simeq \mathrm{Hom}_{S}(M,\mathrm{Hom}%
_{T}(U,N))\text{ for all }M\in \text{ }_{S}\mathbb{M}\text{ and }N\in \text{
}_{T}\mathbb{M}
\end{equation*}%
and the natural transformations%
\begin{equation*}
\omega _{U}^{l}:U\otimes _{S}\mathrm{Hom}_{T}(U,-)\rightarrow 1_{_{T}\mathbb{%
M}}\text{ and }\eta _{U}^{l}:1_{_{S}\mathbb{M}}\rightarrow \mathrm{Hom}%
_{T}(U,U\otimes _{S}-)
\end{equation*}%
yield for every $_{T}K$ and $_{S}L$ the canonical morphisms
\begin{equation}
\omega _{U,K}^{l}:U\otimes _{S}\mathrm{Hom}_{T}(U,K)\rightarrow K\text{ and }%
\eta _{U,L}^{l}:L\rightarrow \mathrm{Hom}_{T}(U,U\otimes _{S}L).
\label{v-eta}
\end{equation}%
We call the $\mathbf{H}_{U}^{l}$\emph{-static }modules $U$\textbf{-static
w.r.t. }$S$ and set
\begin{equation*}
\mathrm{Stat}^{l}(_{T}U_{S}):=\mathrm{Stat}(\mathbf{H}_{U}^{l})=\{_{T}K\mid
U\otimes _{S}\mathrm{Hom}_{T-}(U,K)\overset{\omega _{U,K}^{l}}{\simeq }K\};
\end{equation*}%
and the $\mathbf{H}_{U}^{l}$\emph{-adstatic} modules $U$\textbf{-adstatic
w.r.t. }$S$ and set
\begin{equation*}
\mathrm{Adstat}^{l}(_{T}U_{S}):=\mathrm{Adstat}(\mathbf{H}%
_{U}^{l})=\{_{S}L\mid L\overset{\eta _{U,L}^{l}}{\simeq }\mathrm{Hom}%
_{T-}(U,U\otimes _{S}L)\}.
\end{equation*}%
By \cite{Nau-1990a} and \cite{Nau1990b}, there are equivalences of categories%
\begin{equation}
\mathrm{Stat}^{l}(_{T}U_{S})\approx \mathrm{Adstat}^{l}(_{T}U_{S}).
\label{Ref=Ref}
\end{equation}%
On the other hand, one can define the full subcategories $\mathrm{Stat}%
^{r}(_{T}U_{S})\approx \mathrm{Adstat}^{r}(_{T}U_{S}):$%
\begin{equation*}
\begin{tabular}{lll}
$\mathrm{Stat}^{r}(_{T}U_{S})$ & $:=$ & $\{K_{S}\mid \mathrm{Hom}%
_{-S}(U,K)\otimes _{T}U\simeq K\};$ \\
$\mathrm{Adstat}^{r}(_{T}U_{S})$ & $:=$ & $\{L_{T}\mid L\simeq \mathrm{Hom}%
_{-S}(U,L\otimes _{T}U)\}.$%
\end{tabular}%
\end{equation*}%
In particular, setting%
\begin{equation*}
\begin{tabular}{lllllll}
$\mathrm{Stat}(_{T}U)$ & $:=$ & $\mathrm{Stat}^{l}(_{T}U_{\mathrm{End}%
(_{T}U)^{op}});$ &  & $\mathrm{Adstat}(_{T}U)$ & $:=$ & $\mathrm{Adstat}%
^{l}(_{T}U_{\mathrm{End}(_{T}U)^{op}});$ \\
$\mathrm{Stat}(U_{S})$ & $:=$ & $\mathrm{Stat}^{r}(_{\mathrm{End}%
(_{S}U)}U_{S});$ &  & $\mathrm{Adstat}(U_{S})$ & $:=$ & $\mathrm{Adstat}%
^{r}(_{\mathrm{End}(_{S}U)}U_{S}),$%
\end{tabular}%
\end{equation*}%
there are equivalences of categories:%
\begin{equation}
\mathrm{Stat}(_{T}U)\simeq \mathrm{Adstat}(_{T}U)\text{ and }\mathrm{Stat}%
(U_{S})\simeq \mathrm{Adstat}(U_{S}).  \label{stat=ad}
\end{equation}
\end{punto}

\begin{remark}
The theory of static and adstatic modules was developed in a series of
papers by the second author (see the references). They were also considered
by several other authors (e.g. \cite{Alp1990}, \cite{CF-2004}). For other
terminologies used by different authors, the interested reader may refer to
a comprehensive treatment of the subject by R. Wisbauer in \cite{Wis2000}.
\end{remark}

\subsection*{Intersecting subcategories}

\qquad \qquad Several intersecting subcategories related to Morita contexts
were introduced in the literature (e.g. \cite{Nau-1993}, \cite{Nau-1994-b}).
In what follows we introduce more and we show that many of these coincide,
if one starts with an injective Morita semi-context. Moreover, other results
on equivalences between some intersecting subcategories related to an
injective Morita context will be reframed for arbitrary (not necessarily
compatible) injective Morita data.

\begin{definition}
\begin{enumerate}
\item For a right $T$-module $X,$ a $T$-submodule $X^{\prime }\subseteq X$
is called $K$\textbf{-pure} for some left $T$-module $_{T}K,$ iff the
following sequence of Abelian groups is exact
\begin{equation*}
0\rightarrow X^{\prime }\otimes _{T}K\rightarrow X\otimes _{T}K\rightarrow
X/X^{\prime }\otimes _{T}K\rightarrow 0;
\end{equation*}

\item For a left $T$-module $Y,$ a $T$-submodule $Y^{\prime }\subseteq Y$ is
called $L$\textbf{-copure} for some left $T$-module $_{T}L,$ iff the
following sequence of Abelian groups is exact%
\begin{equation*}
0\rightarrow \mathrm{Hom}_{T}(Y/Y^{\prime },L)\rightarrow \mathrm{Hom}%
_{T}(Y,L)\rightarrow \mathrm{Hom}_{T}(Y^{\prime },L)\rightarrow 0.
\end{equation*}
\end{enumerate}
\end{definition}

\begin{definition}
(Compare \cite[Theorems 1.3., 2.3.]{KO-1979}) Let $T$ be a ring, $%
I\vartriangleleft T$ an ideal, $U$ a left $T$-module and consider the
canonical $T$-linear morphisms%
\begin{equation*}
\zeta _{I,U}:U\rightarrow \mathrm{Hom}_{T}(I,U)\text{ and }\xi
_{I,U}:I\otimes _{T}U\rightarrow U.
\end{equation*}

\begin{enumerate}
\item We say $_{T}U$ is $I$\textbf{-divisible}, iff $\xi _{I,U}$ is
surjective (equivalently, iff $IU=U$).

\item We say $_{T}U$ is $I$\textbf{-localized}, iff $U\overset{\zeta _{I,U}}{%
\simeq }\mathrm{Hom}_{T}(I,U)$ canonically (equivalently iff $_{T}I$ is
strongly $U$-faithful and $_{T}I\subseteq T$ is $U$-copure).

\item We say a left $T$-module $U$ is $I$\textbf{-colocalized}, iff $%
I\otimes _{T}U\overset{\xi _{I,U}}{\simeq }U$ canonically (equivalently, iff
$_{T}U$ is $I$-divisible and $I_{T}\subseteq T$ is $U$-pure).
\end{enumerate}
\end{definition}

\begin{notation}
For a ring $T,$ an ideal $I\vartriangleleft T,$ and with morphisms being the
canonical ones, we set%
\begin{equation*}
\begin{tabular}{lllllll}
$_{I}\mathfrak{D}$ & $:=$ & $\{_{T}U\mid IU=U\};$ &  & $_{I}\mathfrak{F}$ & $%
:=$ & $\{_{T}U\mid U\hookrightarrow \mathrm{Hom}_{T-}(I,U)\};$ \\
$_{I}\mathcal{L}$ & $:=$ & $\{_{T}U\mid U\simeq \mathrm{Hom}_{T}(I,U\};$ &
& $_{I}\mathcal{C}$ & $:=$ & $\{_{T}U\mid I\otimes _{T}U\simeq U\};$ \\
$\mathfrak{D}_{I}$ & $:=$ & $\{U_{T}\mid UI=U\};$ &  & $\mathfrak{F}_{I}$ & $%
:=$ & $\{U_{T}\mid U\hookrightarrow \mathrm{Hom}_{-T}(I,U)\};$ \\
$\mathcal{L}_{I}$ & $:=$ & $\{U_{T}\mid U\simeq \mathrm{Hom}_{T}(I,U\};$ &
& $\mathcal{C}_{I}$ & $:=$ & $\{U_{T}\mid U\otimes _{T}I\simeq U\};.$%
\end{tabular}%
\end{equation*}
\end{notation}

\qquad The following result is due to T. Kato, K. Ohtake and B. M\"{u}ller
(e.g. \cite{Mul-1974}, \cite{Kat-1978}, \cite{KO-1979}):

\begin{proposition}
\label{C=C}\emph{\ }Let $\mathcal{M}=(T,S,P,Q,<,>_{T},<,>_{S},I,J)\in
\mathbb{UMC}.$ Then there are equivalences of categories%
\begin{equation*}
_{I}\mathcal{C}\approx \text{ }_{J}\mathcal{C},\text{ }\mathcal{C}%
_{I}\approx \mathcal{C}_{J},\text{ }_{I}\mathcal{L}\approx \text{ }_{J}%
\mathcal{L}\text{ and }\mathcal{L}_{I}\approx \mathcal{L}_{J}.
\end{equation*}
\end{proposition}

\begin{punto}
\label{UVWX}Let $\mathbf{m}_{T}=(T,S,P,Q,<,>_{T},I)\in \mathbb{UMSC}$ and
consider the dual pairings $\mathbf{P}_{l}:=(Q,$ $_{T}P)\in \mathcal{P}%
_{l}(T)$ and $\mathbf{Q}_{r}:=(P,Q_{T})\in \mathcal{P}_{r}(T).$ For every
left (right) $T$-module $U$ consider the canonical $S$-linear morphism
induced by $<,>_{T}:$%
\begin{equation*}
\alpha _{U}^{\mathbf{Q}_{r}}:Q\otimes _{T}U\rightarrow \mathrm{Hom}_{T-}(P,U)%
\text{ (}\alpha _{U}^{\mathbf{P}_{l}}:U\otimes _{T}P\rightarrow \mathrm{Hom}%
_{-T}(Q,U)\text{).}
\end{equation*}%
We define%
\begin{equation*}
\begin{tabular}{lll}
$\mathcal{D}_{l}(\mathbf{m}_{T})$ & $:=$ & $\{_{T}U\mid Q\otimes _{T}U%
\overset{\alpha _{U}^{\mathbf{Q}_{r}}}{\simeq }\mathrm{Hom}_{T-}(P,U)\};$ \\
$\mathcal{D}_{r}(\mathbf{m}_{T})$ & $:=$ & $\{U_{T}\mid U\otimes _{T}P%
\overset{\alpha _{U}^{\mathbf{P}_{l}}}{\simeq }\mathrm{Hom}_{-T}(Q,U)\}.$%
\end{tabular}%
\end{equation*}%
Moreover, set%
\begin{equation}
\begin{tabular}{lllllll}
$\mathcal{U}_{l}(\mathbf{m}_{T})$ & $:=$ & $\mathrm{Stat}^{l}(_{T}P_{S})\cap
\mathrm{Adstat}^{l}(_{S}Q_{T});$ &  & $\mathcal{U}_{r}(\mathbf{m}_{T})$ & $%
:= $ & $\mathrm{Stat}^{r}(_{S}Q_{T})\cap \mathrm{Adstat}^{r}(_{T}P_{S});$ \\
&  &  &  &  &  &  \\
\QTR{cal}{V}$_{l}(\mathbf{m}_{T})$ & $:=$ & $\mathrm{Stat}%
^{l}(_{T}P_{S})\cap \mathcal{D}_{l}(\mathbf{m}_{T});$ &  & \QTR{cal}{V}$_{r}(%
\mathbf{m}_{T})$ & $:=$ & $\mathrm{Stat}^{r}(_{S}Q_{T})\cap \mathcal{D}_{r}(%
\mathbf{m}_{T});$ \\
$\mathbb{V}_{l}(\mathbf{m}_{T})$ & $:=$ & $_{I}\mathcal{C}\cap \mathcal{D}%
_{l}(\mathbf{m}_{T});$ &  & $\mathbb{V}_{r}(\mathbf{m}_{T})$ & $:=$ & $%
\mathcal{C}_{I}\cap \mathcal{D}_{r}(\mathbf{m}_{T});$ \\
$\widehat{\mathcal{V}}_{l}(\mathbf{m}_{T})$ & $:=$ & $\mathcal{V}_{l}(%
\mathbf{m}_{T})\cap $ $_{I}\mathcal{L};$ &  & $\widehat{\mathcal{V}}_{r}(%
\mathbf{m}_{T})$ & $:=$ & $\mathcal{V}_{r}(\mathbf{m}_{T})\cap \mathcal{L}%
_{I};$ \\
&  &  &  &  & $:=$ &  \\
$\mathcal{W}_{l}(\mathbf{m}_{T})$ & $:=$ & $\mathrm{Adstat}%
^{l}(_{S}Q_{T})\cap \mathcal{D}_{l}(\mathbf{m}_{T});$ &  & $\mathcal{W}_{r}(%
\mathbf{m}_{T})$ & $:=$ & $\mathrm{Adstat}^{r}(_{T}P_{S})\cap \mathcal{D}%
_{r}(\mathbf{m}_{T});$ \\
$\mathbb{W}_{l}(\mathbf{m}_{T})$ & $:=$ & $_{I}\mathcal{L}\cap \mathcal{D}%
_{l}(\mathbf{m}_{T});$ &  & $\mathbb{W}_{r}(\mathbf{m}_{T})$ & $:=$ & $%
\mathcal{L}_{I}\cap \mathcal{D}_{r}(\mathbf{m}_{T});$ \\
$\widehat{\mathcal{W}}_{l}(\mathbf{m}_{T})$ & $:=$ & $\mathcal{W}_{l}(%
\mathbf{m}_{T})\cap $ $_{I}\mathcal{C};$ &  & $\widehat{\mathcal{W}}_{r}(%
\mathbf{m}_{T})$ & $:=$ & $\mathcal{W}_{r}(\mathbf{m}_{T})\cap \mathcal{C}%
_{I};$ \\
&  &  &  &  &  &  \\
$\mathcal{X}_{l}(\mathbf{m}_{T})$ & $:=$ & $\mathcal{V}_{l}(\mathbf{m}%
_{T})\cap \mathcal{W}_{l}(\mathbf{m}_{T});$ &  & $\mathcal{X}_{r}(\mathbf{m}%
_{T})$ & $:=$ & $\mathcal{V}_{r}(\mathbf{m}_{T})\cap \mathcal{W}_{r}(\mathbf{%
m}_{T});$ \\
$\mathbb{X}_{l}(\mathbf{m}_{T})$ & $:=$ & $\mathbb{V}_{l}(\mathbf{m}%
_{T})\cap \mathbb{W}_{l}(\mathbf{m}_{T});$ &  & $\mathbb{X}_{r}(\mathbf{m}%
_{T})$ & $:=$ & $\mathbb{V}_{r}(\mathbf{m}_{T})\cap \mathbb{W}_{r}(\mathbf{m}%
_{T}).$ \\
&  &  &  &  &  &  \\
$\mathcal{X}_{l}^{\ast }(\mathbf{m}_{T})$ & $:=$ & $\{_{S}(Q\otimes
_{T}U)\mid V\in \mathcal{X}_{l}(\mathbf{m}_{T})\};$ &  & $\mathcal{X}%
_{r}^{\ast }(\mathbf{m}_{T})$ & $:=$ & $\{(U\otimes _{T}P)_{S}\mid V\in
\mathcal{X}_{r}(\mathbf{m}_{T})\};$ \\
$\mathbb{X}_{l}^{\ast }(\mathbf{m}_{T})$ & $:=$ & $\{_{S}(Q\otimes
_{T}U)\mid V\in \mathbb{X}_{l}(\mathbf{m}_{T})\};$ &  & $\mathbb{X}%
_{r}^{\ast }(\mathbf{m}_{T})$ & $:=$ & $\{(U\otimes _{T}P)_{S}\mid V\in
\mathbb{X}_{r}(\mathbf{m}_{T})\}.$%
\end{tabular}
\label{LIST}
\end{equation}%
Given $\mathbf{m}_{S}=(S,T,Q,P,<,>_{S},J)\in \mathbb{UMSC}$ one can define
analogously, the corresponding intersecting subcategories of $_{S}\mathbb{M}$
and $\mathbb{M}_{S}.$
\end{punto}

As an immediate consequence of Proposition \ref{C=C} we get

\begin{corollary}
\label{V=C}Let $\mathcal{M}=(T,S,P,Q,<,>_{T},<,>_{S},I,J)\in \mathbb{IUMC}$
and consider the associated Morita semi-contexts $\mathcal{M}_{T}$ and $%
\mathcal{M}_{S}$ \emph{(\ref{MT-MS})}.

\begin{enumerate}
\item If $_{I}\mathcal{C}\leq \mathcal{D}_{l}(\mathcal{M}_{T})$ and $_{J}%
\mathcal{C}\leq \mathcal{D}_{l}(\mathcal{M}_{S}),$ then $\mathcal{V}_{l}(%
\mathcal{M}_{T})\approx \mathcal{V}_{l}(\mathcal{M}_{S}).$ Similarly, if $%
\mathcal{C}_{I}\leq \mathcal{D}_{r}(\mathcal{M}_{T})$ and $\mathcal{C}%
_{J}\leq \mathcal{D}_{r}(\mathcal{M}_{S}),$ then $\mathcal{V}_{r}(\mathcal{M}%
_{T})\approx \mathcal{V}_{r}(\mathcal{M}_{S}).$

\item If $_{I}\mathcal{L}\leq \mathcal{D}_{l}(\mathcal{M}_{T})$ and $_{J}%
\mathcal{L}\leq \mathcal{D}_{l}(\mathcal{M}_{S}),$ then $\mathcal{W}_{l}(%
\mathcal{M}_{T})\approx \mathcal{W}_{l}(\mathcal{M}_{S}).$ Similarly, if $%
\mathcal{L}_{I}\leq \mathcal{D}_{r}(\mathcal{M}_{T})$ and $\mathcal{L}%
_{J}\leq \mathcal{D}_{r}(\mathcal{M}_{S}),$ then $\mathcal{W}_{r}(\mathcal{M}%
_{T})\approx \mathcal{W}_{r}(\mathcal{M}_{S}).$
\end{enumerate}
\end{corollary}

\qquad Starting with a Morita context, the following result was obtained in
\cite[Theorem 3.2.]{Nau-1993}. We restate the result for an arbitrary (not
necessarily compatible) Morita datum and \emph{sketch} its proof:

\begin{lemma}
\label{X=X}Let $\mathcal{M}=\mathbf{(}T,S,P,Q,<,>_{T},<,>_{S},I,J)$ be a
unital Morita datum and consider the associated Morita semi-contexts $%
\mathcal{M}_{T}$ and $\mathcal{M}_{S}$ in \emph{(\ref{MT-MS})}. Then there
are equivalences of categories%
\begin{equation*}
\mathcal{X}_{l}(\mathcal{M}_{T})\overset{\mathrm{Hom}_{T-}(P,-)}{\underset{%
\mathrm{Hom}_{S-}(Q,-)}{\approx }}\mathcal{X}_{l}(\mathcal{M}_{S})\text{ and
}\mathcal{X}_{r}(\mathcal{M}_{T})\overset{\mathrm{Hom}_{-T}(Q,-)}{\underset{%
\mathrm{Hom}_{-S}(P,-)}{\approx }}\mathcal{X}_{r}(\mathcal{M}_{S}).
\end{equation*}
\end{lemma}

\begin{Beweis}
Let $_{T}V\in \mathcal{X}_{l}(\mathcal{M}_{T}).$ By the equivalence $\mathrm{%
Stat}^{l}(_{T}P_{S})\overset{\mathrm{Hom}_{T}(P,-)}{\approx }\mathrm{A}%
\mathrm{dstat}^{l}(_{T}P_{S})$ in \ref{ref=ref} we have $\mathrm{Hom}%
_{T-}(P,V)\in \mathrm{Adstat}^{l}(_{T}P_{S}).$ Moreover, $V\in \mathcal{D}%
_{l}(M),$ hence $\mathrm{Hom}_{T-}(P,V)\simeq Q\otimes _{T}V$ canonically
and it follows then from the equivalence $\mathrm{Adstat}^{l}(_{S}Q_{T})%
\overset{Q\otimes _{T}-}{\approx }\mathrm{Stat}^{l}(_{S}Q_{T})$ that $%
\mathrm{Hom}_{T-}(P,V)\in \mathrm{Stat}^{l}(_{S}Q_{T}).$ Moreover, we have
the following \emph{natural} isomorphisms%
\begin{equation}
P\otimes _{S}\mathrm{Hom}_{T-}(P,V)\simeq V\simeq \mathrm{Hom}%
_{S-}(Q,Q\otimes _{T}V)\simeq \mathrm{Hom}_{S-}(Q,\mathrm{Hom}_{T-}(P,V)),
\label{Hom-D}
\end{equation}%
i.e. $\mathrm{Hom}_{T-}(P,V)\in \mathcal{D}_{l}(\mathcal{M}_{S}).$
Consequently, $\mathrm{Hom}_{T-}(P,V)\in \mathcal{X}_{l}(\mathcal{M}_{S}).$
Moreover, (\ref{Hom-D}) yields a natural isomorphism $V\simeq \mathrm{Hom}%
_{S-}(Q,\mathrm{Hom}_{T-}(P,V)).$ Analogously, one can show for every $W\in
\mathcal{X}_{l}(\mathcal{M}_{S})$ that $\mathrm{Hom}_{S-}(Q,W)\in \mathcal{X}%
_{l}(\mathcal{M}_{T})$ and that $W\simeq \mathrm{Hom}_{T-}(P,\mathrm{Hom}%
_{S-}(Q,W))$ naturally. Consequently, $\mathcal{X}_{l}(\mathcal{M}%
_{T})\approx \mathcal{X}_{l}(\mathcal{M}_{S})$. The equivalences $\mathcal{X}%
_{r}(\mathcal{M}_{T})\approx \mathcal{X}_{r}(\mathcal{M}_{S})$ can be proved
analogously.$\blacksquare $
\end{Beweis}

\begin{proposition}
\label{Stat-I-I}Let $\mathcal{M}=(T,S,P,Q,<,>_{T},<,>_{S},I,J)$ be a unital
injective Morita datum and consider the associated Morita semi-contexts $%
\mathcal{M}_{T}$ and $\mathcal{M}_{S}$ in \emph{(\ref{MT-MS}).}

\begin{enumerate}
\item There are equivalences of categories%
\begin{equation*}
\begin{tabular}{lllllll}
$\mathrm{Stat}^{l}(_{T}I_{T})$ & $\approx $ & $\mathrm{Adstat}%
^{l}(_{T}I_{T});$ &  & $\mathrm{Stat}^{l}(_{S}J_{S})$ & $\approx $ & $%
\mathrm{Adstat}^{l}(_{S}J_{S});$ \\
$\mathrm{Stat}^{r}(_{T}I_{T})$ & $\approx $ & $\mathrm{Adstat}%
^{r}(_{T}I_{T});$ &  & $\mathrm{Stat}^{r}(_{S}J_{S})$ & $\approx $ & $%
\mathrm{Adstat}^{r}(_{S}J_{S}).$%
\end{tabular}%
\end{equation*}

\item If $\mathrm{Stat}^{l}(_{T}I_{T})\leq \mathcal{X}_{l}^{\ast }(\mathcal{M%
}_{S})$ and $\mathrm{Stat}^{l}(_{S}J_{S})\leq \mathcal{X}_{l}^{\ast }(%
\mathcal{M}_{T}),$ then there are equivalences of categories%
\begin{equation*}
\mathrm{Stat}^{l}(_{T}I_{T})\approx \mathrm{Stat}^{l}(_{S}J_{S})\text{ and }%
\mathrm{Adstat}^{l}(_{T}I_{T})\approx \mathrm{Adstat}^{l}(_{S}J_{S}).
\end{equation*}

\item If $\mathrm{Stat}^{r}(_{T}I_{T})\leq \mathcal{X}_{r}^{\ast }(\mathcal{M%
}_{S})$ and $\mathrm{Stat}^{r}(_{S}J_{S})\leq \mathcal{X}_{r}^{\ast }(%
\mathcal{M}_{T}),$ then there are equivalences of categories%
\begin{equation*}
\mathrm{Stat}^{r}(_{T}I_{T})\approx \mathrm{Stat}^{r}(_{S}J_{S})\text{ and }%
\mathrm{Adstat}^{r}(_{T}I_{T})\approx \mathrm{Adstat}^{r}(_{S}J_{S}).
\end{equation*}
\end{enumerate}
\end{proposition}

\begin{Beweis}
To prove \textquotedblleft 1\textquotedblright , notice that since $\mathcal{%
M}$ is an injective Morita datum, $P\otimes _{S}Q\overset{<,>_{T}}{\simeq }I$
and $Q\otimes _{T}P\overset{<,>_{S}}{\simeq }J$ as bimodules and so the four
equivalences of categories result from \ref{ref=ref}. To prove
\textquotedblleft 2\textquotedblright , one can use an argument similar to
that in \cite[Theorem 3.9.]{Nau-1994-b} to show that the inclusion $\mathrm{%
Stat}^{l}(_{T}I_{T})=\mathrm{Stat}^{l}(_{T}(P\otimes _{S}Q)_{T})\leq
\mathcal{X}_{l}^{\ast }(\mathcal{M}_{S})$ implies $\mathrm{Stat}%
^{l}(_{T}I_{T})=\mathrm{Stat}^{l}(_{T}(P\otimes _{S}Q)_{T})=\mathcal{X}_{l}(%
\mathcal{M}_{T})$ and that the inclusion $\mathrm{Stat}^{l}(_{S}J_{S})=%
\mathrm{Stat}^{l}(_{S}(Q\otimes _{T}P)_{S})\leq \mathcal{X}_{l}^{\ast }(%
\mathcal{M}_{T})$ implies $\mathrm{Stat}^{l}(_{S}J_{S})=\mathrm{Stat}%
^{l}(_{S}(Q\otimes _{T}P)_{S})=\mathcal{X}_{l}(\mathcal{M}_{S}).$ The result
follows then by Lemma \ref{X=X}. The proof of \textquotedblleft
3\textquotedblright\ is analogous to that of \textquotedblleft
2\textquotedblright .$\blacksquare $
\end{Beweis}

\qquad For injective Morita semi-contexts, several subcategories in (\ref%
{LIST}) are shown in the following result to be equal:

\begin{theorem}
\label{V=V}Let $\mathbf{m}_{T}=(T,S,P,Q,<,>_{T},I)\in \mathbb{IUMS}.$ Then

\begin{enumerate}
\item $\mathcal{V}_{l}(\mathbf{m}_{T})=\mathbb{V}_{l}(\mathbf{m}_{T}),$ $%
\mathcal{W}_{l}(\mathbf{m}_{T})=\mathbb{W}_{l}(\mathbf{m}_{T}),$ whence
\begin{equation*}
\widehat{\mathcal{V}}_{l}(\mathbf{m}_{T})=\widehat{\mathcal{W}}_{l}(\mathbf{m%
}_{T})=\mathcal{X}_{l}(\mathbf{m}_{T})=\mathbb{X}_{l}(\mathbf{m}_{T})=\text{
}_{I}\mathcal{C}\cap \mathcal{D}_{l}(\mathbf{m}_{T})\cap \text{ }_{I}%
\mathcal{L}\text{ and }\mathcal{X}_{l}^{\ast }(\mathbf{m}_{T})=\mathbb{X}%
_{l}^{\ast }(\mathbf{m}_{T}).
\end{equation*}

\item $\mathcal{V}_{r}(\mathbf{m}_{T})=\mathbb{V}_{r}(\mathbf{m}_{T}),$ $%
\mathcal{W}_{r}(\mathbf{m}_{T})=\mathbb{W}_{r}(\mathbf{m}_{T}),$ whence
\begin{equation*}
\widehat{\mathcal{V}}_{r}(\mathbf{m}_{T})=\widehat{\mathcal{W}}_{r}(\mathbf{m%
}_{T})=\mathcal{X}_{r}(\mathbf{m}_{T})=\mathbb{X}_{r}(\mathbf{m}_{T})=%
\mathcal{C}_{I}\cap \mathcal{D}_{r}(\mathbf{m}_{T})\cap \mathcal{L}_{I}\text{
and }\mathcal{X}_{r}^{\ast }(\mathbf{m}_{T})=\mathbb{X}_{r}^{\ast }(\mathbf{m%
}_{T}).
\end{equation*}
\end{enumerate}
\end{theorem}

\begin{Beweis}
We prove only \textquotedblleft 1\textquotedblright\ as \textquotedblleft
2\textquotedblright\ can be proved analogously. Assume the Morita
semi-context $\mathbf{m}_{T}=(T,S,P,Q,<,>_{T},I)$ is injective. By our
assumption we have for every $V\in \mathcal{D}_{l}(\mathbf{m}_{T})$ the
commutative diagram%
\begin{equation}
\xymatrix{ P\otimes _{S}(Q\otimes _{T}V) \ar[rrrr]^{can} _{\simeq}
\ar[dd]_{id_{P}\otimes _{S}(\alpha _{V}^{\mathbf{Q}_{r}})}^{\simeq} & & & &
(P\otimes _{S}Q)\otimes _{T}V \ar[dd]^{{<,>_{T}\otimes _{T}id_{V}}}_{\simeq}
\\ & & & &\\ P\otimes _{S}\mathrm{Hom}_{T-}(P,V) \ar[rr]_(.7){\omega
_{P,V}^{l}} & & V & & I\otimes _{_{T}}V \ar[ll]^{\xi _{I,V}}}  \label{S}
\end{equation}%
Then it becomes obvious that $\omega _{P,V}^{l}:P\otimes _{S}\mathrm{Hom}%
_{T}(P,V)\rightarrow V$ is an isomorphism if and only if $\xi
_{I,V}:I\otimes _{T}V\rightarrow V$ is an isomorphism. Consequently%
\begin{equation*}
\mathcal{V}(\mathbf{m}_{T})=\mathcal{D}_{l}(\mathbf{m}_{T})\cap \mathrm{Stat}%
^{l}(_{T}P_{S})=\mathcal{D}_{l}(\mathbf{m}_{T})\cap \text{ }_{I}\mathcal{C}=%
\mathbb{V}(\mathbf{m}_{T}).
\end{equation*}%
On the other hand, we have for every $V\in \mathcal{D}_{l}(\mathbf{m}_{T})$
the following commutative diagram%
\begin{equation}
\xymatrix{ \mathrm{Hom}_{S-}(Q,\mathrm{Hom}_{T-}(P,V)) \ar[rrrr]^{can}
_{\simeq} & & & & \mathrm{Hom}_{T-}(P\otimes _{S}Q,V) \\ & & & & \\
\mathrm{Hom}_{S-}(Q,Q\otimes _{T}V) \ar[uu]^{(Q,\alpha
_{V}^{\mathbf{Q}_{r}})}_{\simeq} & & V \ar[ll]^(.3){\eta_{P,L}^{l}}
\ar[rr]_{\zeta_{I,V}}& & \mathrm{Hom}_{T-}(I,V)
\ar[uu]_{(<,>_{T},V)}^{\simeq} }  \label{Ad}
\end{equation}%
It follows then that $\eta _{P,L}^{l}:V\rightarrow \mathrm{Hom}%
_{S}(Q,Q\otimes _{T}P)$ is an isomorphism if and only if $\zeta
_{I,V}:V\rightarrow \mathrm{Hom}_{T}(I,V)$ is an isomorphism. Consequently,
\begin{equation*}
\mathcal{W}(\mathbf{m}_{T})=\mathcal{D}_{l}(\mathbf{m}_{T})\cap \mathrm{%
Adstat}^{l}(_{T}P_{S})=\mathcal{D}_{l}(\mathbf{m}_{T})\cap _{I}\mathcal{L}=%
\mathbb{W}(\mathbf{m}_{T}).
\end{equation*}%
Moreover, we have%
\begin{equation*}
\begin{tabular}{lllllll}
$\widehat{\mathcal{V}}_{l}(\mathbf{m}_{T})$ & $:=$ & $\mathcal{V}_{l}(%
\mathbf{m}_{T})\cap \text{ }_{I}\mathcal{L}$ & $=$ & $\mathbb{V}_{l}(\mathbf{%
m}_{T})\cap \text{ }_{I}\mathcal{L}$ & $=$ & $_{I}\mathcal{C}\cap \mathcal{D}%
_{l}(\mathbf{m}_{T})\cap $ $_{I}\mathcal{L}$ \\
& $=$ & $_{I}\mathcal{C}\cap \mathbb{W}_{l}(\mathbf{m}_{T})$ & $=$ & $_{I}%
\mathcal{C}\cap \mathcal{W}_{l}(\mathbf{m}_{T})$ & $=$ & $\widehat{\mathcal{W%
}}_{l}(\mathbf{m}_{T}).$%
\end{tabular}%
\end{equation*}%
On the other hand, we have%
\begin{equation*}
\mathcal{X}_{l}(\mathbf{m}_{T})=\mathcal{V}_{l}(\mathbf{m}_{T})\cap \mathcal{%
W}_{l}(\mathbf{m}_{T})=\mathbb{V}_{l}(\mathbf{m}_{T})\cap \mathbb{W}_{l}(%
\mathbf{m}_{T})=\mathbb{X}_{l}(\mathbf{m}_{T})
\end{equation*}%
and so the equalities $\widehat{\mathcal{V}}_{l}(\mathbf{m}_{T})=\widehat{%
\mathcal{W}}_{l}(\mathbf{m}_{T})=\mathcal{X}_{l}(\mathbf{m}_{T})=\mathbb{X}%
_{l}(\mathbf{m}_{T})$ and $\mathcal{X}_{l}^{\ast }(\mathbf{m}_{T})=\mathbb{X}%
_{l}^{\ast }(\mathbf{m}_{T})$ are established.$\blacksquare $
\end{Beweis}

\qquad In addition to establishing several other equivalences of
intersecting subcategories, the following results reframe the equivalence of
categories $\widehat{\mathcal{V}}\approx \widehat{\mathcal{W}}$ in \cite[%
Theorem 4.9.]{Nau-1994-b} for an arbitrary (not necessarily compatible)
injective Morita datum:

\begin{theorem}
\label{CHECK}Let $\mathcal{M}=(T,S,P,Q,<,>_{T},<,>_{S},I,J)$ be an injective
Morita datum and consider the associated Morita semi-contexts $\mathcal{M}%
_{T}$ and $\mathcal{M}_{S}$ \emph{(\ref{MT-MS})}.

\begin{enumerate}
\item The following subcategories are mutually equivalent:
\begin{equation}
\widehat{\mathcal{V}}_{l}(\mathcal{M}_{T})=\widehat{\mathcal{W}}_{l}(%
\mathcal{M}_{T})=\mathbb{X}_{l}(\mathcal{M}_{T})=\mathcal{X}_{l}(\mathcal{M}%
_{T})\approx \mathcal{X}_{l}(\mathcal{M}_{S})=\mathbb{X}_{l}(\mathcal{M}%
_{S})=\widehat{\mathcal{W}}_{l}(\mathcal{M}_{S})=\widehat{\mathcal{V}}_{l}(%
\mathcal{M}_{S}).  \label{8T}
\end{equation}

\item If $\mathcal{V}_{l}(\mathcal{M}_{T})\leq $ $_{I}\mathcal{L}$ and $%
\mathcal{W}_{l}(\mathcal{M}_{S})\leq \mathcal{\ }_{J}\mathcal{C},$ then $%
\mathcal{V}_{l}(\mathcal{M}_{T})\approx \mathcal{W}_{l}(\mathcal{M}_{S}).$
If $\mathcal{W}_{l}(\mathcal{M}_{T})\leq $ $_{I}\mathcal{C}$ and $\mathcal{V}%
_{l}(\mathcal{M}_{S})\leq \mathcal{\ }_{J}\mathcal{L},$ then $\mathcal{W}%
_{l}(\mathcal{M}_{T})\approx \mathcal{V}_{l}(\mathcal{M}_{S}).$

\item The following subcategories are mutually equivalent:
\begin{equation}
\widehat{\mathcal{V}}_{r}(\mathcal{M}_{T})=\widehat{\mathcal{W}}_{r}(%
\mathcal{M}_{T})=\mathbb{X}_{r}(\mathcal{M}_{T})=\mathcal{X}_{r}(\mathcal{M}%
_{T})\approx \mathcal{X}_{r}(\mathcal{M}_{S})=\mathbb{X}_{r}(\mathcal{M}%
_{S})=\widehat{\mathcal{W}}_{r}(\mathcal{M}_{S})=\widehat{\mathcal{V}}_{r}(%
\mathcal{M}_{S}).
\end{equation}

\item If $\mathcal{V}_{r}(\mathcal{M}_{T})\leq \mathcal{L}_{I}$ and $%
\mathcal{W}_{r}(\mathcal{M}_{T})\leq \mathcal{C}_{J},$ then $\mathcal{V}_{r}(%
\mathcal{M}_{T})\approx \mathcal{W}_{r}(\mathcal{M}_{S}).$ If $\mathcal{W}%
_{r}(\mathcal{M}_{T})\leq \mathcal{C}_{J}$ and $\mathcal{V}_{r}(\mathcal{M}%
_{S})\leq \mathcal{L}_{I},$ then $\mathcal{V}_{r}(\mathcal{M}_{S})\approx
\mathcal{W}_{r}(\mathcal{M}_{T}).$
\end{enumerate}
\end{theorem}

\begin{Beweis}
By Lemma \ref{X=X}, $\mathcal{X}_{l}(\mathcal{M}_{T})\approx \mathcal{X}_{l}(%
\mathcal{M}_{S})$ and so \textquotedblleft 1\textquotedblright\ follows by
Theorem \ref{V=V}. If $\mathcal{V}_{l}(\mathcal{M}_{T})\leq $ $_{I}\mathcal{L%
}$ and $\mathcal{W}_{l}(\mathcal{M}_{S})\leq \mathcal{\ }_{J}\mathcal{C},$
then we have%
\begin{equation*}
\mathcal{V}_{l}(\mathcal{M}_{T})=\mathcal{V}_{l}(\mathcal{M}_{T})\cap \text{
}_{I}\mathcal{L}=\widehat{\mathcal{V}}_{l}(\mathcal{M}_{T})\approx \widehat{%
\mathcal{W}}_{l}(\mathcal{M}_{S})=\mathcal{W}_{l}(\mathcal{M}_{S})\cap \text{
}_{J}\mathcal{C}=\mathcal{W}_{l}(\mathcal{M}_{S}).
\end{equation*}%
On the other hand, if $\mathcal{W}_{l}(\mathcal{M}_{T})\leq $ $_{I}\mathcal{L%
}$ and $\mathcal{V}_{l}(\mathcal{M}_{S})\leq \mathcal{\ }_{J}\mathcal{C},$
then%
\begin{equation*}
\mathcal{W}_{l}(\mathcal{M}_{T})=\mathcal{W}_{l}(\mathcal{M}_{T})\cap \text{
}_{I}\mathcal{C}=\widehat{\mathcal{W}}_{l}(\mathcal{M}_{T})\approx \widehat{%
\mathcal{V}}_{l}(\mathcal{M}_{S})=\mathcal{V}_{l}(\mathcal{M}_{S})\cap \text{
}_{J}\mathcal{L}=\mathcal{V}_{l}(\mathcal{M}_{S}).
\end{equation*}%
So we have established \textquotedblleft 2\textquotedblright . The results
in \textquotedblleft 3\textquotedblright\ and \textquotedblleft
4\textquotedblright\ can be obtained analogously.$\blacksquare $
\end{Beweis}

\section{More applications}

\qquad \qquad In this final section we give more applications of Morita $%
\alpha $-(semi-)contexts and injective Morita (semi-)contexts. All rings in
this section are \emph{unital}, whence all Morita (semi-)contexts are
unital. Moreover, for any ring $T$ we denote with $_{T}\mathbf{E}$ an
arbitrary, but fixed, injective cogenerator in $_{T}\mathbb{M}.$

\begin{notation}
Let $T$ be an $A$-ring. For any left $T$-module $_{T}V,$ we set $^{\#}V:=%
\mathrm{Hom}_{T}(V,$ $_{T}\mathbf{E}).$ If moreover, $_{T}V_{S}$ is a $(T,S)$%
-bimodule for some $B$-ring $S,$ then we consider $_{S}^{\#}V$\emph{\ }with
the left $S$-module structure induced by that of $V_{S}.$
\end{notation}

\begin{lemma}
\label{gen-surj}\emph{(Compare \cite[Lemma 3.2.]{Col1990}, \cite[Lemmas
2.1.2., 2.1.3.]{CF-2004})}\ Let $T$ be an $A$-ring, $S$ a $B$-ring and $%
_{T}V_{S}$ a $(T,S)$-bimodule,

\begin{enumerate}
\item A left $T$-module $_{T}K$ is $V$-generated if and only if the
canonical $T$-linear morphism
\begin{equation}
\omega _{V,K}^{l}:V\otimes _{S}\mathrm{Hom}_{T}(V,K)\rightarrow K
\label{omega}
\end{equation}%
is surjective. Moreover, $V\otimes _{S}W\subseteq \mathrm{Pres}%
(_{T}V)\subseteq \mathrm{Gen}(_{T}V)$ for every left $S$-module $_{S}W.$

\item A left $S$-module $_{S}L$ is $_{S}^{\#}V$-cogenerated if and only if
the canonical $S$-linear morphism
\begin{equation}
\eta _{V,L}^{l}:L\rightarrow \mathrm{Hom}_{T}(V,V\otimes _{S}L)  \label{eta}
\end{equation}%
is injective. Moreover, $\mathrm{Hom}_{T}(V,M)\subseteq \mathrm{Copres}%
(_{S}^{\#}V)\subseteq \mathrm{Cogen}(_{S}^{\#}V)$ for every left $T$-module $%
_{T}M.$
\end{enumerate}
\end{lemma}

\begin{remark}
Let $T$ be an $A$-ring, $S$ a $B$-ring and $_{T}V_{S}$ a $(T,S)$-bimodule.
Notice that for any left $S$-module $_{S}L$ we have%
\begin{equation*}
\mathrm{ann}_{L}^{\otimes }(V_{S}):=\{l\in L\mid V\otimes _{S}l=0\}=\mathrm{%
Ker}(\eta _{V,L}^{l}),
\end{equation*}%
whence (by Lemma \ref{gen-surj} \textquotedblleft 2\textquotedblright\ ) $%
V_{S}$ is $L$-faithful if and only if $_{S}L$ is $_{S}^{\#}V$-cogenerated.
It follows then that $V_{S}$ is completely faithful if and only if $%
_{S}^{\#}V$ is a cogenerator.
\end{remark}

\subsection*{Localization and colocalization}

\qquad In what follows we clarify the relations between static (adstatic)
modules and subcategories colocalized (localized) by a trace ideal of a
Morita context satisfying the $\alpha $-condition.

Recall that for any $(T,S)$-bimodule $_{T}P_{S}$ we have by Lemma \ref%
{gen-surj}:%
\begin{equation}
\mathrm{Stat}^{l}(_{T}P_{S})\subseteq \mathrm{Gen}(_{T}P)\text{ and }\mathrm{%
Adstat}^{l}(_{T}P_{S})\subseteq \mathrm{Cogen}(_{S}^{\#}P).  \label{stat-gen}
\end{equation}

\begin{theorem}
\label{Gen=reflex}Let $\mathcal{M}=(T,S,P,Q,<,>_{T},<,>_{S},I,J)\in \mathbb{%
UMC}.$ Then we have%
\begin{equation}
_{I}\mathcal{C}\subseteq \text{ }_{I}\mathfrak{D}\subseteq \mathrm{Gen}%
(_{T}P).  \label{REF=REF}
\end{equation}%
Assume $\mathbf{P}_{r}:=(Q,P_{S})\in \mathcal{P}_{r}^{\alpha }(S).$ Then

\begin{enumerate}
\item $\mathrm{Gen}(_{T}P)=\mathrm{Stat}^{l}(_{T}P_{S})\subseteq $ $_{I}%
\mathfrak{F}.$

\item If $\mathrm{Gen}(_{T}P)\subseteq $ $_{I}\mathcal{C},$ then $_{I}%
\mathcal{C}=$ $_{I}\mathfrak{D}=\mathrm{Gen}(_{T}P)=\mathrm{Stat}%
^{l}(_{T}P_{S}).$

\item If $\mathbf{Q}_{r}:=(P,Q_{T})\in \mathcal{P}_{r}^{\alpha }(T),$ then $%
_{T}I\subseteq $ $_{T}T$ is pure and $_{I}\mathcal{C}=$ $_{I}\mathfrak{D}.$
\end{enumerate}
\end{theorem}

\begin{Beweis}
For every left $T$-module $_{T}K,$ consider the following diagram with
canonical morphisms and let $\alpha _{2}:=\zeta _{I,K}\circ \omega
_{P,K}^{l}.$ It is easy to see that both rectangles and the two right
triangles commutes:%
\begin{equation}
{\xymatrix{P \otimes_S Q \otimes_T K \ar^{id_P \otimes _S \alpha_K^{{\bold
Q}_r}}[rr] \ar_{<,>_{T} \otimes_T id_K}[dd] & & P \otimes_S {\rm Hom}_T(P,K)
\ar^{\alpha^{{\bold P}_{r}}_{{\rm Hom}_T(P,K)}}[rr]
\ar^{\omega_{P,K}^{l}}[dd] \ar@{.>}[rrdd]_{\alpha_2} & & {\rm
Hom}_{S}(Q,{\rm Hom}_T(P,K)) \\& & & & {\rm Hom}_T(P \otimes_S Q,K)
\ar_{\simeq}[u] \\ I \otimes_{T} K \ar_{\xi_{I,K}}[rr]
\ar@{.>}[rruu]_{\alpha_1} & & K \ar_{\zeta_{I,K}}[rr] & & {\rm Hom}_T(I,K)
\ar_{(<,>_T,K)}[u]}}  \label{I}
\end{equation}%
It follows directly from the definitions that $_{I}\mathcal{C}\subseteq $ $%
_{I}\mathfrak{D}$ and $\mathrm{Stat}^{l}(_{T}P_{S})\subseteq \mathrm{Gen}%
(_{T}P).$ If $_{T}K$ is $I$-divisible, then $\xi _{I,K}\circ <,>_{T}\otimes
_{T}id_{K}=\omega _{P,K}^{l}\circ id_{P}\otimes _{S}\alpha _{K}^{\mathbf{Q}%
_{r}}$ is surjective, whence $\omega _{P,K}^{l}$ is surjective and we
conclude that $_{T}K$ is $P$-generated by Lemma \ref{gen-surj}
\textquotedblleft 1\textquotedblright . Consequently, $_{I}\mathfrak{D}%
\subseteq \mathrm{Gen}(_{T}P).$

Assume now that $\mathbf{P}_{r}\in \mathcal{P}_{r}^{\alpha }(S).$
Considering the canonical map $\rho _{Q}:T\rightarrow \mathrm{End}%
(_{S}Q)^{op},$ the map $\rho _{Q}\circ <,>_{T}=\alpha _{Q}^{\mathbf{P}_{r}}$
is injective and so the bilinear map $<,>_{T}$ is injective (i.e. $P\otimes
_{S}Q\overset{<,>_{T}}{\simeq }I$). Define $\alpha _{1}:=(id_{P}\otimes
_{S}\alpha _{K}^{\mathbf{Q}_{r}})\circ (<,>_{T}\otimes _{T}id_{K})^{-1},$ so
that the left triangles commute. Notice that $\alpha _{\mathrm{Hom}%
_{T}(P,K)}^{\mathbf{P}_{r}}$ is injective and the commutativity of the upper
right triangle in Diagram (\ref{I}) implies that $\alpha _{2}$ is injective
(whence $\omega _{P,K}^{l}$ is injective by the commutativity of the lower
right triangle).

\begin{enumerate}
\item If $K\in \mathrm{Stat}^{l}(_{T}P_{S}),$ then the commutativity of the
lower right triangle (\ref{I})\ and the injectivity of $\alpha _{2}$ show
that $\zeta _{I,K}$ is injective; hence, $\mathrm{Stat}^{l}(_{T}P_{S})%
\subseteq $ $_{I}\mathfrak{F}.$ On the other hand, if $_{T}K$ is $P$%
-generated, then $\omega _{P,K}^{l}$ is surjective by Lemma \ref{gen-surj}
(1), thence bijective, i.e. $K\in \mathrm{Stat}^{l}(_{T}P_{S}).$
Consequently, $\mathrm{Gen}(_{T}P)=\mathrm{Stat}^{l}(_{T}P_{S}).$

\item This follows directly from the inclusions in (\ref{REF=REF}) and
\textquotedblleft 1\textquotedblright .

\item Assume $\mathbf{Q}_{r}:=(P,Q_{T})\in \mathcal{P}_{r}^{\alpha }(T).$
Since $\mathbf{P}_{r}\in \mathcal{P}_{r}^{\alpha }(S),$ it follows by
analogy to Proposition \ref{rp-rp} \textquotedblleft 3\textquotedblright\
that $P_{S}$ is flat, hence $id_{P}\otimes _{S}\alpha _{K}^{\mathbf{Q}_{r}}$
is injective. The commutativity of the upper left triangle in Diagram (\ref%
{I}) implies then that $\alpha _{1}$ is injective, thence $\xi _{I,K}$ is
injective by commutativity of the lower left triangle (i.e. $_{T}I\subseteq $
$_{T}T$ is $K$-pure). If $_{T}K$ is divisible, then $K\otimes _{T}I\overset{%
\xi _{I,K}}{\simeq }K$ (i.e. $K\in $ $_{I}\mathcal{C}$).$\blacksquare $
\end{enumerate}
\end{Beweis}

\begin{theorem}
\label{Cog=ref}Let $\mathcal{M}=(T,S,P,Q,<,>_{T},<,>_{S},I,J)\in \mathbb{UMC}%
.$ Then we have%
\begin{equation*}
_{J}\mathcal{L}\subseteq \text{ }_{J}\mathfrak{F}\subseteq \mathrm{Cogen}%
(_{S}^{\#}P)\text{ and }\mathrm{Adstat}^{l}(_{T}P_{S})\subseteq \mathrm{Cogen%
}(_{S}^{\#}P).
\end{equation*}%
Assume $\mathbf{Q}_{r}:=(P,Q_{T})\in \mathcal{P}_{r}^{\alpha }(T).$ Then

\begin{enumerate}
\item $J_{S}\subseteq S_{S}$ is pure and $_{J}\mathcal{C}\subseteq \mathrm{%
Cogen}(_{S}^{\#}P).$

\item If $\mathbf{P}_{r}:=(Q,P_{S})\in \mathcal{P}_{r}^{\alpha }(S),$ then $%
_{J}\mathcal{L}\subseteq \mathrm{Adstat}^{l}(_{T}P_{S})\subseteq \mathrm{%
Cogen}(_{S}^{\#}P)\subseteq $ $_{J}\mathfrak{F}.$

\item If $\mathbf{P}_{r}\in \mathcal{P}_{r}^{\alpha }(S)$ and $\mathrm{Cogen}%
(_{S}^{\#}P)\subseteq $ $_{J}\mathcal{L},$ then $_{J}\mathcal{L}=\mathrm{%
Cogen}(_{S}^{\#}P)=\mathrm{Adstat}^{l}(_{T}P_{S}).$
\end{enumerate}
\end{theorem}

\begin{Beweis}
For every right $S$-module $L$ consider the commutative diagram with
canonical morphisms and let $\alpha _{3}$ be so defined, that the left
triangles become commutative%
\begin{equation}
\xymatrix{J \otimes_S L\ar^{\xi_{J,L}}[rr] \ar@{.>}[rrdd]_{\alpha_3} & &
L\ar_{\eta_{P,L}^{l}}[dd] \ar^{\zeta_{J,L}}[rr] & & {\rm Hom}_S(J,L)
\ar^{(<,>_S,L)}[d]\\ & & & & {\rm Hom}_S(Q \otimes_T P,L) \ar^{\simeq \, \,
{\rm can}}[d]\\ Q \otimes_T P \otimes_S L \ar^{(<,>_{S}) \otimes_S id_L
}[uu] \ar_{\alpha^{{\bold Q}_r}_{P \otimes_S L}}[rr] & & {\rm Hom}_T (P,
P\otimes_S L) \ar[rr]_{(P,\alpha^{{\bold P}_r}_L)} \ar@{.>}[rruu]_{\alpha_4}
& & {\rm Hom}_T(P,{\rm Hom}_S(Q,L))}  \label{J}
\end{equation}%
By definition $_{J}\mathcal{L}\subseteq $ $_{J}\mathfrak{F}$ and $\mathrm{%
Adstat}^{l}(_{T}P_{S})\subseteq \mathrm{Cogen}(_{S}^{\#}P).$ If $_{S}L\in $ $%
_{J}\mathfrak{F},$ then $\zeta _{J,L}$ is injective and it follows by
commutativity of the right rectangle in Diagram (\ref{J}) that $\eta
_{P,L}^{l}$ is injective, hence $_{S}L$ is $_{S}^{\#}P$-cogenerated by Lemma %
\ref{gen-surj} \textquotedblleft 2\textquotedblright . Consequently, $_{J}%
\mathfrak{F}\subseteq \mathrm{Cogen}(_{S}^{\#}P).$

Assume now that $\mathbf{Q}_{r}\in \mathcal{P}_{r}^{\alpha }(T).$ Then it
follows from Lemma \ref{alpha->inj} that $<,>_{S}$ is injective (hence $%
Q\otimes _{T}P\overset{<,>_{S}}{\simeq }J$) and so $\alpha _{4}:=(\mathrm{can%
}\circ (<,>_{S},L))^{-1}\circ (P,\alpha _{L}^{\mathbf{P}_{r}})$ is injective.

\begin{enumerate}
\item Since $\alpha _{3}$ is injective, $\xi _{J,L}$ is also injective for
every $_{S}L,$ i.e. $J_{S}\subseteq S_{S}$ is pure. If $_{S}L\in $ $_{J}%
\mathcal{C},$ then it follows from the commutativity of the left rectangle
in Diagram (\ref{J}) that $\eta _{P,L}^{l}$ is injective, hence $L\in
\mathrm{Cogen}(_{S}^{\#}P)$ by Lemma \ref{gen-surj} (2).

\item Assume that $\mathbf{P}_{r}\in \mathcal{P}_{r}^{\alpha }(S),$ so that $%
\alpha _{4}$ is injective. If $_{S}L\in $ $_{J}\mathcal{L},$ then $\zeta
_{J,L}$ is an isomorphism, thence $\eta _{P,L}^{l}$ is surjective (notice
that $\alpha _{4}$ is injective). Consequently, $_{J}\mathcal{L}\subseteq
\mathrm{Adstat}^{l}(_{T}P_{S}).$

\item This follows directly from the assumptions and \textquotedblleft
2\textquotedblright .$\blacksquare $
\end{enumerate}
\end{Beweis}

\subsection*{$\ast $-Modules}

\qquad To the end of this section, we fix a unital ring $T,$ a left $T$%
-module $_{T}P$ and set $S:=\mathrm{End}(_{T}P)^{op}.$

\bigskip

\begin{definition}
(\cite{MO-1989}) We call $_{T}P$ a $\ast $\textbf{-module,} iff $\mathrm{Gen}%
(_{T}P)\approx \mathrm{Cogen}(_{S}^{\#}P).$
\end{definition}

\begin{remark}
It was shown by J. Trlifaj \cite{Trl1994} that all $\ast $-modules are
finitely generated.
\end{remark}

\qquad By definition, $\mathrm{Stat}^{l}(_{T}P_{S})\leq $ $_{T}\mathbb{M}$
and $\mathrm{Adstat}^{l}(_{T}P_{S})\leq $ $_{S}\mathbb{M}$ are the \emph{%
largest} subcategories between which the adjunction $(P\otimes _{S}-,\mathrm{%
Hom}_{T}(P,-))$ induces an equivalence. On the other hand, Lemma \ref%
{gen-surj} shows that $\mathrm{Gen}(_{T}P)\leq $ $_{T}\mathbb{M}$ and $%
\mathrm{Cogen}(_{S}^{\#}P)\leq $ $_{S}\mathbb{M}$ are the \emph{largest}
such subcategories (see \cite[Section 3]{Col1990} for more details). This
suggests the following observation:

\begin{proposition}
\label{star-stat}\emph{(\cite[Lemma 2.3.]{Xin1999})} We have\emph{\ }%
\begin{equation*}
_{T}P\text{ is a }\ast \text{-module}\Leftrightarrow \mathrm{Stat}(_{T}P)=%
\mathrm{Gen}(_{T}P)\text{ and }\mathrm{Adstat}(_{T}P)=\mathrm{Cogen}%
(_{S}^{\#}P).
\end{equation*}
\end{proposition}

\begin{definition}
A left $T$-module $_{T}U$ is said to be

\textbf{semi-}$\dsum $\textbf{-quasi-projective }(abbr. $s$\textbf{-}$\dsum $%
\textbf{-quasi-projective}), iff for any left $T$-module $_{T}V\in \mathrm{%
Pres}(_{T}U)$ and any $U$\emph{-presentation}%
\begin{equation*}
U^{(\Lambda )}\rightarrow U^{(\Lambda ^{\prime })}\rightarrow V\rightarrow 0
\end{equation*}%
of $_{T}V$(if any), the following induced sequence is exact:%
\begin{equation*}
\mathrm{Hom}_{T}(U,U^{(\Lambda )})\rightarrow \mathrm{Hom}_{T}(U,U^{(\Lambda
^{\prime })})\rightarrow \mathrm{Hom}_{T}(U,V)\rightarrow 0;
\end{equation*}

\textbf{weakly-}$\dsum $\textbf{-quasi-projective }(abbr. $w$\textbf{-}$%
\dsum $\textbf{-quasi-projective}), iff for any left $T$-module $_{T}V$ and
any short exact sequence%
\begin{equation*}
0\rightarrow K\rightarrow U^{(\Lambda ^{\prime })}\rightarrow V\rightarrow 0
\end{equation*}%
with $K\in \mathrm{Gen}(_{T}U)$ (if any), the following induced sequence is
exact:%
\begin{equation*}
0\rightarrow \mathrm{Hom}_{T}(U,K)\rightarrow \mathrm{Hom}_{T}(U,U^{(\Lambda
^{\prime })})\rightarrow \mathrm{Hom}_{T}(U,V)\rightarrow 0;
\end{equation*}

\textbf{self-tilting}, iff $_{T}U$ is $w$-$\dsum $-quasi-projective and $%
\mathrm{Gen}(_{T}U)=\mathrm{Pres}(_{T}U);$

$\dsum $\textbf{-self-static}, iff any direct sum $U^{(\Lambda )}$ is $U$%
-static.

(\textbf{self})\textbf{-small, }iff $\mathrm{Hom}_{T}(U,-)$ commutes with
direct sums (of $_{T}U$);
\end{definition}

\begin{proposition}
\label{gen=stat}Assume $\mathcal{M}=(T,S,P,Q,<,>_{T},<,>_{S})$ is a unital
Morita context.

\begin{enumerate}
\item If $\mathbf{P}_{r}:=(Q,P_{S})\in \mathcal{P}_{r}^{\alpha }(S),$ then:

\begin{enumerate}
\item $\mathrm{Gen}(_{T}P)=\mathrm{Stat}^{l}(_{T}P_{S});$

\item there is an equivalence of categories $\mathrm{Gen}(_{T}P)\approx
\mathrm{Cop}(_{S}^{\#}P);$

\item $_{T}P$ is $\dsum $-self-static and $\mathrm{Stat}^{l}(_{T}P_{S})$ is
closed under factor modules.

\item $\mathrm{Gen}(_{T}P)=\mathrm{Pres}(_{T}P);$
\end{enumerate}

\item If $\mathcal{M}\in \mathbb{UMC}_{r}^{\alpha }$ and $\mathrm{Cogen}%
(_{S}^{\#}P)\subseteq $ $_{J}\mathcal{L},$ then:

\begin{enumerate}
\item $\mathrm{Gen}(_{T}P)=\mathrm{Stat}^{l}(_{T}P_{S})$ and $\mathrm{Cogen}%
(_{S}^{\#}P)=\mathrm{Adstat}^{l}(_{T}P_{S});$

\item there is an equivalence of categories $\mathrm{Cogen}%
(_{S}^{\#}P)\approx \mathrm{Gen}(_{T}P);$

\item $_{T}P$ is a $\ast $-module;

\item $_{T}P$ is self-tilting and self-small.
\end{enumerate}
\end{enumerate}
\end{proposition}

\begin{Beweis}
\begin{enumerate}
\item If $\mathbf{P}_{r}\in \mathcal{P}_{r}^{\alpha }(S),$ then it follows
by Theorem \ref{Gen=reflex} that $\mathrm{Gen}(_{T}P)=\mathrm{Stat}%
^{l}(_{T}P_{S}),$ which is equivalent to each of \textquotedblleft
b\textquotedblright\ and \textquotedblleft c\textquotedblright\ by \cite[4.4.%
]{Wis2000} and to \textquotedblleft d\textquotedblright\ by \cite[4.3.]%
{Wis2000}.

\item It follows by the assumptions, Theorems \ref{Gen=reflex}, \ref{Cog=ref}
and \ref{ref=ref} that $\mathrm{Gen}(_{T}P)=\mathrm{Stat}^{l}(_{T}P_{S})%
\approx \mathrm{Adstat}^{l}(_{T}P_{S})=\mathrm{Cogen}(_{S}^{\#}P),$ whence $%
\mathrm{Gen}(_{T}P)\approx \mathrm{Cogen}(_{S}^{\#}P)$ (which is the
definition of $\ast $-modules). Hence \textquotedblleft a\textquotedblright\
$\Leftrightarrow $\textquotedblleft b\textquotedblright\ $\Leftrightarrow $%
\textquotedblleft c\textquotedblright . The equivalence \textquotedblleft
a\textquotedblright\ $\Leftrightarrow $ \textquotedblleft
d\textquotedblright\ is evident by \cite[Corollary 4.7.]{Wis2000} and we are
done.$\blacksquare $
\end{enumerate}
\end{Beweis}

\subsection*{Wide Morita Contexts}

\qquad \emph{Wide Morita contexts} were introduced by F. Casta\~{n}o
Iglesias and J. G\'{o}mez-Torrecillas \cite{C-IG-T1995} and \cite{C-IG-T1996}
as an extension of classical \emph{Morita contexts} to Abelian categories.

\begin{definition}
Let $\mathcal{A}$ and $\mathcal{B}$ be Abelian categories. A \textbf{right }(%
\textbf{left})\textbf{\ wide Morita context} between $\mathcal{A}$ and $%
\mathcal{B}$ is a datum $\mathcal{W}_{r}=(G,\mathcal{A},\mathcal{B},F,\eta
,\rho ),$ where $G:\mathcal{A}\rightleftarrows \mathcal{B}:F$ are right
(left) exact covariant functors and $\eta :F\circ G\longrightarrow 1_{%
\mathcal{A}},$ $\rho :G\circ F\longrightarrow 1_{\mathcal{B}}$ ($\eta :1_{%
\mathcal{A}}\longrightarrow F\circ G,$ $\rho :1_{\mathcal{B}}\longrightarrow
G\circ F$) are natural transformations, such that for every pair of objects $%
(A,B)\in \mathcal{A}\times \mathcal{B}$ the compatibility conditions $G(\eta
_{A})=\rho _{G(A)}$ and $F(\rho _{B})=\eta _{F(B)}\ $hold.
\end{definition}

\begin{definition}
Let $\mathcal{A}$ and $\mathcal{B}$ be Abelian categories and $\mathcal{W}%
=(G,\mathcal{A},\mathcal{B},F,\eta ,\rho )$ be a right (left) wide Morita
context. We call $\mathcal{W}$ \textbf{injective }(respectively \textbf{%
semi-strict}, \textbf{strict}), iff $\eta $ and $\rho $ are monomorphisms
(respectively epimorphisms, isomorphisms)
\end{definition}

\begin{remarks}
Let $\mathcal{W}=(G,\mathcal{A},\mathcal{B},F,\eta ,\rho )$ be a right
(left) wide Morita context.

\begin{enumerate}
\item It follows by \cite[Propositions 1.1., 1.4.]{CDN2005} that if either $%
\eta $ or $\rho $ is an epimorphism (monomorphism), then $\mathcal{W}$ is
strict, whence $\mathcal{A}\approx \mathcal{B}.$

\item The resemblance of \emph{injective} left wide Morita contexts\textit{\
}is with the Morita-Takeuchi contexts for comodules of coalgebras, i.e. the
so called \emph{pre-equivalence data} for categories of comodules introduced
in \cite{Tak1977} (see \cite{C-IG-T1998} for more details).
\end{enumerate}
\end{remarks}

\subsubsection*{Injective Right wide Morita contexts}

\qquad

In a recent work \cite[5.1.]{CDN2005}, Chifan, et. al. clarified (for module
categories) the relation between \emph{classical Morita contexts} and \emph{%
right wide Morita contexts}\textit{. }For the convenience of the reader and
for later reference, we include in what follows a brief description of this
relation.

\begin{punto}
Let $T,S$ be rings, $\mathcal{A}:=$ $_{T}\mathbb{M}$ and $\mathcal{B}:=$ $%
_{S}\mathbb{M}.\ $Associated to each Morita context\textit{\ }$\mathcal{M}%
=(T,S,P,Q,<,>_{T},<,>_{S})$ is a wide Morita context as follows: Define $G:%
\mathcal{A}\rightleftarrows \mathcal{B}:F$ by $G(-)=Q\otimes _{T}-$ and $%
F(-)=P\otimes _{S}-.$ Then there are natural transformations $\eta :F\circ
G\longrightarrow 1_{_{T}\mathbb{M}}$ and $\rho :G\circ F\longrightarrow
1_{_{S}\mathbb{M}}$ such that for each $_{T}V$ and $W_{S}:$%
\begin{equation}
\begin{tabular}{llllllll}
$\eta _{V}$ & $:$ & $P\otimes _{S}(Q\otimes _{T}V)$ & $\rightarrow $ & $V,$
& $\dsum p_{i}\otimes _{S}(q_{i}\otimes _{T}v_{i})$ & $\mapsto $ & $\dsum
<p_{i},q_{i}>_{T}v_{i},$ \\
$\rho _{W}$ & $:$ & $Q\otimes _{T}(P\otimes _{S}W)$ & $\rightarrow $ & $W,$
& $\dsum q_{i}\otimes _{T}(p_{i}\otimes _{S}w_{i})$ & $\mapsto $ & $\dsum
<q_{i},p_{i}>_{S}w_{i}.$%
\end{tabular}
\label{3.1.}
\end{equation}%
Then the datum $\mathcal{W}_{r}(\mathcal{M}):=(G,$ $_{T}\mathbb{M},$ $_{S}%
\mathbb{M},F,\eta ,\rho )$ is a right wide Morita context\textit{.}

Conversely, let $T^{\prime },S^{\prime }$ be two rings and $\mathcal{W}%
_{r}^{\prime }=(G^{\prime },$ $_{T^{\prime }}\mathbb{M},$ $_{S^{\prime }}%
\mathbb{M},F^{\prime },\eta ^{\prime },\rho ^{\prime })$\textit{\ be a right
wide }Morita context\textit{\ between }$_{T^{\prime }}\mathbb{M}$ and $%
_{S^{\prime }}\mathbb{M}$ such that \textit{the right exact functors }$%
G^{\prime }:$ $_{T^{\prime }}\mathbb{M}\rightleftarrows $ $_{S^{\prime }}%
\mathbb{M}:F^{\prime }$ \textit{commute with direct sums. By \emph{Watts'
Theorems} \emph{(e.g. \cite{Gol1979})},\ }there exists a\textit{\ }$(T,S)$%
-bimodule $P^{\prime }$ (e.g. $F^{\prime }(S^{\prime }))$ such that $%
F^{\prime }\simeq P^{\prime }\otimes _{S^{\prime }}-,$ an $(S,T)$-bimodule $%
Q^{\prime }$ such that $G^{\prime }\simeq Q^{\prime }\otimes _{T^{\prime }}-$
and there should exist two bilinear forms%
\begin{equation*}
<,>_{T^{\prime }}:P^{\prime }\otimes _{S^{\prime }}Q^{\prime }\rightarrow
T^{\prime }\text{ and }<,>_{S^{\prime }}:Q^{\prime }\otimes _{T^{\prime
}}P^{\prime }\rightarrow S^{\prime },
\end{equation*}%
such that the natural transformations $\eta ^{\prime }:F^{\prime }\circ
G^{\prime }\rightarrow 1_{_{T^{\prime }}\mathbb{M}},$ $\rho :G^{\prime
}\circ F^{\prime }\rightarrow 1_{_{S^{\prime }}\mathbb{M}}$ are given by%
\begin{equation*}
\eta _{V^{\prime }}^{\prime }(p^{\prime }\otimes _{S^{\prime }}q^{\prime
}\otimes _{T^{\prime }}v^{\prime })=<p^{\prime },q^{\prime }>_{T^{\prime
}}v^{\prime }\text{ and }\rho _{W^{\prime }}^{\prime }(q^{\prime }\otimes
_{T}p^{\prime }\otimes _{S}w^{\prime })=<q^{\prime },p^{\prime }>_{S^{\prime
}}w^{\prime }
\end{equation*}%
for all $V^{\prime }\in $ $_{T^{\prime }}\mathbb{M},$ $W^{\prime }\in $ $%
_{S^{\prime }}\mathbb{M},$ $p^{\prime }\in P^{\prime },$ $q^{\prime }\in
Q^{\prime },$ $v^{\prime }\in V^{\prime }$ and $w^{\prime }\in W^{\prime }.$
It can be shown that in this way one obtains a Morita context $\mathcal{M}%
^{\prime }=\mathcal{M}^{\prime }(\mathcal{W}_{r}^{\prime }):=(T^{\prime
},S^{\prime },P^{\prime },Q^{\prime },<,>_{T^{\prime }},<,>_{S^{\prime }}).$
\textit{Moreover, it turns out that given a wide Morita context }$\mathcal{W}%
_{r},$ we have $\mathcal{W}_{r}\simeq \mathcal{W}_{r}(\mathcal{M}(\mathcal{W}%
_{r})).$
\end{punto}

\qquad The following result clarifies the relation between \emph{injective
Morita contexts} and \emph{injective right wide Morita contexts.}

\begin{theorem}
\label{MW-inj}Let $\mathcal{M}=(T,S,P,Q,<,>_{T},<,>_{S})$ be a Morita
context, $\mathcal{A}:=$ $_{T}\mathbb{M},$ $\mathcal{B}:=$ $_{S}\mathbb{M}$
and consider the induced right wide Morita context $\mathcal{W}_{r}(\mathcal{%
M}):=(G,\mathcal{A},\mathcal{B},F,\eta ,\rho ).$

\begin{enumerate}
\item \textit{If }$\mathcal{W}_{r}(\mathcal{M})$ \textit{is an }injective
right wide Morita context\textit{, then }$\mathcal{M}$ \textit{is an }%
injective Morita context\textit{.}

\item If $\mathcal{M}\in \mathbb{UMC}_{r}^{\alpha },$ then\textit{\ }$%
\mathcal{W}_{r}(\mathcal{M})$ \textit{is an }injective right wide Morita
context\textit{.}
\end{enumerate}
\end{theorem}

\begin{Beweis}
\begin{enumerate}
\item Let $\mathcal{W}_{r}(\mathcal{M})$ be an injective right wide Morita
context.\textit{\ }Then in particular, $<,>_{T}=\eta _{T}$ and $<,>_{S}=\rho
_{S}$ are injective, i.e. $\mathcal{M}$ is an injective Morita context.

\item Assume that $\mathcal{M}$ satisfies the right $\alpha $-condition.
Suppose there exists some $_{T}V$ and $\dsum p_{i}\otimes _{S}(q_{i}\otimes
_{T}v_{i})\in \mathrm{Ker}(\eta _{V}).$ Then for any $q\in Q$ we have%
\begin{equation*}
\text{%
\begin{tabular}{lllll}
$0$ & $=$ & $q\otimes _{T}\eta _{V}(\dsum (p_{i}\otimes _{S}q_{i})\otimes
_{T}v_{i})$ & $=$ & $\dsum q\otimes _{T}<p_{i},q_{i}>_{T}v_{i}$ \\
& $=$ & $\dsum q<p_{i},q_{i}>_{T}\otimes _{T}v_{i}$ & $=$ & $\dsum
<q,p_{i}>_{S}q_{i}\otimes _{T}v_{i}$ \\
& $=$ & $\dsum <q,p_{i}>_{S}(q_{i}\otimes _{T}v_{i})$ & $=$ & $\alpha
_{Q\otimes _{T}V}^{\mathbf{P}_{r}}(\dsum p_{i}\otimes _{S}(q_{i}\otimes
_{T}v_{i}))(q).$%
\end{tabular}%
}
\end{equation*}%
Since $\mathbf{P}_{r}:=(Q,P_{S})\in \mathcal{P}_{r}^{\alpha }(S),$ the
morphism $\alpha _{Q\otimes _{T}V}^{\mathbf{P}_{r}}$ is injective and so $%
\dsum p_{i}\otimes _{S}(q_{i}\otimes _{T}v_{i})=0,$ i.e. $\eta _{V}$ is
injective. Analogously, suppose $\dsum q_{i}\otimes _{T}(p_{i}\otimes
_{S}w_{i})\in \mathrm{Ker}(\rho _{W}).$ Then for any $p\in P$ we have%
\begin{equation*}
\begin{tabular}{lllll}
$0$ & $=$ & $p\otimes _{S}\rho _{W}(\dsum q_{i}\otimes _{T}(p_{i}\otimes
_{S}w_{i})$ & $=$ & $\dsum p\otimes _{S}<q_{i},p_{i}>_{S}w_{i}$ \\
& $=$ & $\dsum p<q_{i},p_{i}>_{S}\otimes _{S}w_{i}$ & $=$ & $\dsum
<p,q_{i}>_{T}p_{i}\otimes _{S}w_{i}$ \\
& $=$ & $\dsum <p,q_{i}>_{T}(p_{i}\otimes _{S}w_{i})$ & $=$ & $\alpha
_{P\otimes _{S}W}^{\mathbf{Q}_{r}}(\dsum q_{i}\otimes _{T}(p_{i}\otimes
_{S}w_{i}))(p).$%
\end{tabular}%
\end{equation*}%
Since $\mathbf{Q}_{r}:=(P,Q_{T})\in \mathcal{P}_{r}^{\alpha }(T),$ the
morphism $\alpha _{P\otimes _{S}W}^{\mathbf{Q}_{r}}$ is injective and so $%
\dsum q_{i}\otimes _{T}(p_{i}\otimes _{S}w_{i})=0,$ i.e. $\rho _{W}$ is
injective. Consequently, the induced right wide Morita context $\mathcal{W}%
_{r}(\mathcal{M})$ is injective.$\blacksquare $
\end{enumerate}
\end{Beweis}

\textbf{Acknowledgement: }The authors thank the referee for his/her careful
reading of the paper and for the fruitful suggestions, comments and
corrections, which helped in improving several parts of the paper. Moreover,
they acknowledge the excellent research facilities as well as the support of
their respective institutions, King Fahd University of Petroleum and
Minerals and King AbdulAziz University.

\end{document}